\documentclass[12pt,A4]{article}
\usepackage[utf8]{inputenc}
\usepackage{xcolor}
\usepackage{amssymb}
\usepackage{hyperref}
\usepackage{amsfonts}
\usepackage{ytableau}
\usepackage{mathtools}  
\usepackage{amsmath,   amsthm, graphicx, 
endnotes ,setspace,
verbatim, geometry, times, helvet, courier, bm, url, babel, dcolumn}

\newcommand{\ppmm}[1]{{#1}}\newcommand{\ppmx}[1]{}

\newcommand{\pmatrixx}[1]{\begin{pmatrix} #1 \end{pmatrix}} 
\newcommand{\C}{{\mathbf{C}}}

\def\b{\beta}
\newcommand{\bse}{\begin{subequations}}
\newcommand{\ese}{\end{subequations}}

\newcommand{\R}{\check{R}}

\newcommand{\beq}{\begin{equation}}
\newcommand{\eeq}{\end{equation}}

\newcommand{\bea}{\begin{eqnarray*}}
\newcommand{\eea}{\end{eqnarray*}}

\newcommand{\ket}[1]{| #1 \rangle}
\newcommand{\bra}[1]{ \langle #1 |}

\newcommand{\ignore}[1]{}

\newtheorem{theorem}{Theorem}[section]
\newtheorem{proposition}[theorem]{Proposition}
\newtheorem{lemma}[theorem]{Lemma}
\newtheorem{corollary}[theorem]{Corollary}
\newtheorem{conjecture}[theorem]{Conjecture}
\newcommand{\propo}[1]{\begin{proposition} #1 \end{proposition}}

\newtheorem{defin}[theorem]{Definition}

\newcommand{\mythsf}{}

\newcommand{\aaaaa}{{\mythsf a}}  \newcommand{\aaaaaa}{{a}}
\newcommand{\bbbbb}{{\mythsf b}}  \newcommand{\bbbbbb}{{b}}

\newcommand{\ccccc}{{\mythsf c}}  \newcommand{\ddddd}{{\mythsf d}}  \newcommand{\dddddd}{{d}}  
\newcommand{\xxxxx}{{\mythsf x}} 
\newcommand{\x}{{\mythsf x}}

\title{Solutions to the constant Yang-Baxter equation:
\\ 
additive charge conservation in three dimensions
}

\author{J.Hietarinta\footnote{
  hietarin@utu.fi, University of Turku, Finland }, 
P.Martin\footnote{
  p.p.martin@leeds.ac.uk, University of Leeds, UK },
 E.C.Rowell\footnote{
rowell@tamu.edu, Texas A\& M University, USA  }}

\begin{document}

\maketitle
\begin{abstract}
We find all solutions to the constant Yang--Baxter equation
$R_{12}R_{13}R_{23}=R_{23}R_{13}R_{12}$ in three dimensions, subject
to an additive charge-conservation ansatz.  This ansatz is a
generalisation of (strict) charge-conservation, for which a complete
classification in all dimensions was recently obtained. Additive
charge-conservation introduces additional sector-coupling parameters
-- in 3 dimensions there are $4$ such parameters.  In the generic
dimension 3 case, in which all of the $4$ parameters are nonzero, we
find there is a single 3 parameter family of solutions.  We give a complete
analysis of this solution, giving the structure of the centraliser
(symmetry) algebra in all orders. We also solve the remaining cases
with three, two, or one nonzero sector-coupling parameter(s).
\end{abstract}

\section{Introduction}

\newcommand{\RR}{R}          \newcommand{\calR}{{\cal R}}

The Yang-Baxter equation (YBE) reads (in shorthand form)
\beq\label{eq:YBs2}
\RR_{12}R_{13}R_{23}=R_{23}R_{13}R_{12} .
\eeq
It is
a fundamental equation for many applications
--- see for example
\cite{Baxter72,Baxter73,STF1979,Baxter1982},
\cite{Drinfeld1985,Jimbo1985,Jones1990,Turaev1988}
and references therein.

To make \eqref{eq:YBs2} explicit, one first fixes
a dimension $N$ for a vector space $V=\C^N$.
\ppmm{We can also pick bases for $V$ and $V \otimes V$.}
Then we have an underlying
matrix
$R$ acting on
$V \otimes V$.
Each matrix $R_{ij}$ acts on $V \otimes V \otimes V$,  
acting on  the $i$-th and $j$-th factors
as $R$,
and on the
other factor as the 
identity.
Thus in
explicit form the
Yang-Baxter equation reads
\begin{equation}\label{eq:YB}
 \sum_{\alpha_1,\alpha_2,\alpha_3}{\calR}^{i_1i_2}_{\alpha_1\alpha_2}
     {\cal R}^{\alpha_1i_3}_{j_1\alpha_3} {\cal
       R}^{\alpha_2\alpha_3}_{j_2j_3} =
     \sum_{\beta_1,\beta_2,\beta_3}{\cal R}^{i_2i_3}_{\beta_2\beta_3}
      {\cal R}^{i_1\beta_3}_{\beta_1j_3} {\cal
       R}^{\beta_1\beta_2}_{j_1j_2},
\end{equation}
where  the indices range over $0,1,\dots,N-1 \;$
and $ \calR^{i_1i_2}_{\alpha_1\alpha_2}$
is the appropriate matrix entry of $R$.
(See also \S\ref{ss:presenting}.)

With various applications in mind, we impose
\beq \det(R) \neq 0.  \eeq

For some applications the $R$ matrices depend on
spectral parameters
that can be different for each $R_{ij}$ \cite{STF1979,Baxter1982},
but in this paper we will  consider the
constant YBE.  By construction, any $\R$ gives a representation of the braid group $B_n$ for each $n$.

Observe that $\RR$ 
will have $N^2\times N^2$ entries and
there will be, in principle, $N^3\times N^3$ equations. It is clear
that such an overdetermined set of nonlinear equations is difficult to
solve, even in this
constant form.
Indeed, while many individual solutions are known,
a complete solution is known
only for dimension two \cite{Hietarinta1992} and for higher dimensions
{  knowledge is far from complete}.
The three dimensional upper triangular case was
solved in \cite{Hietarinta_1993}, but for further progress it is
important to make a meaningful ansatz.

Recently Martin and Rowell proposed \cite{MR} charge-conservation
of the form
\beq
\calR_{ij}^{kl}=0,\text{ if } \{i,j\}\neq\{k,l\} \text{
  as a set,}
\eeq
as an effective constraint and with it they were
able to find all solutions for all dimensions. The above constraint
may be called ``strict charge conservation'' (SCC). In this paper we
will explore the results obtained by relaxing the SCC rule to
``additive charge conservation'' (ACC) defined by
\beq\label{eq:ACC}
\calR_{ij}^{kl}=0,\text{ if } i+j\neq k+l.
\eeq
Observe that 
ACC differs from SCC first
in dimension 3.  In practice this change increases the complexity of
the underlying computational problem by introducing four further
`mixing' parameters
(SCC itself having fifteen parameters in dimension 3).

The paper is organized as follows. In Section \ref{ss:YB} we discuss
notational matters and symmetries of the problem. In Section
\ref{S:sol} we present the solutions.  The  set of solutions is
organized according to the non-vanishing conditions on the four mixing
parameters (together with their symmetries).  Thus the first family
of solutions is the generic case, with all parameters non-zero --
it is solved in detail in \S\ref{ss:Jxxxx}.  The various possibilities
are then addressed in turn, the last case being the set of solutions where
all but one mixing parameter vanishes - \S\ref{ss:case6}.

It turns out that  in several solutions
have the `Hecke' property (i.e., having precisely two distinct eigenvalues).  In
\S\ref{ss:analys} we use this to analyse the representations, giving a
complete analysis for the generic case.

\medskip 

{
One natural realisation of the constant Yang--Baxter problem is as a problem in
categorical representation theory, and this is the perspective largely
taken in \cite{MR} (see also \cite{MRTorzewska}, for example). 
However here we will keep to a simple analytical setting.
Direct transliteration of results between the settings is a routine exercise.
}

\noindent
{\bf{Acknowledgements}}.
We thank Frank Nijhoff for various important contributions,
including initiating our collaboration. 
PM thanks EPSRC for funding under grant EP/W007509/1;
and Paula Martin for useful conversations. ECR was partially funded by US NSF grants DMS-2000331 and DMS-2205962.

\section{The setup}  \label{ss:YB}

\ppmx{
[Re  \eqref{eq:braidr}
  --- if we do not retain that eqn above, which we might not, then we'd need to
adjust wording here, and  copy it here] ---
The $\R$ in braid relation \eqref{eq:braidr} and $R$ in YBE
\eqref{eq:YBs2} are related by}

{For the braid group point of view we first define}
\beq\label{eq:brex}
\R=P\, R,\,\text{ where }\, \mathcal{P}_{ij}^{kl}=\delta_{i}^{l}\delta_{j}^{k},
\quad\text{ i.e. }\, \mathcal{\check{R}}_{ij}^{kl}=\calR_{ji}^{kl}.
\eeq
and furthermore
$$ (P R)_{12} \; = \; \R_1 \; := \; \R \otimes 1
\quad \mbox{ and } \quad (P R)_{23} =\R_2 \; := \; 1\otimes \R
$$
acting on $V\otimes V\otimes V$.
Then the YBE in \eqref{eq:YBs2} becomes
\beq\label{eq:braidr}
\R_1\R_2\R_1=\R_2\R_1\R_2,
\eeq
i.e., the braid group version of the YBE.

\subsection{Presenting matrices}\label{ss:presenting}

Set $V=\C^3$ with basis labeled by $\{0,1,2\}$.
{We will order this basis as the symbols suggest.}
Using the standard ket notation, i.e.
$i\otimes j=:\ket{ij}$, we
may order the basis of $V\otimes V$ for example using lexicographical order
\[
\ket{00}, \ket{01}, \ket{02}, \ket{10}, \ket{11}, \ket{12}, \ket{20}, \ket{21}, \ket{22}
\]
 or reverse lexicographical order (rlex)

\[
\ket{00}, \ket{10}, \ket{20}, \ket{01}, \ket{11}, \ket{21},\ket{ 02}, \ket{12}, \ket{22}
\]
Still another possibility is to use
{ a }
`graded'
reverse lexicographical ordering (grlex)
\[
\ket{00}, \ket{10}, \ket{01}, \ket{20}, \ket{11}, \ket{02}, \ket{21}, \ket{12}, \ket{22}
\]
{The name is borrowed from monomial orderings, in which setting
the symbols {\em are} numbers, rather than being arbitrarily associated to numbers as in our case.}

The matrix entries are defined as:
$$
\calR_{ij}^{kl}:=\bra{ij}R\ket{kl}
$$

In the present case with ACC \eqref{eq:ACC} and the rlex ordering we get the matrix

\begin{equation}\label{eq:matf}
  {R_{{rlex}}  =
\begin{pmatrix}{\cal R}_{0,0}^{0,0} &. & . & .
 & . & . & .  & . & .\\ . & {\cal R}_{1,0}^{1,0} & . & {\cal
    R}_{1,0}^{0,1} & . & . & . & . & .\\ . & . & {\cal R}_{2,0}^{2,0}
  & . & {\cal R}_{2,0}^{1,1} & . & {\cal R}_{2,0}^{0,2 } & . &
  . \\ . & {\cal R}_{0,1}^{1,0} & . & {\cal R}_{0,1}^{0,1} & . & . &
  . & . & . \\ . & . & {\cal R}_{1,1}^{2,0} & . & {\cal R}_{1,1}^{1,1}
  & . & {\cal R}_{1,1}^{0,2} & .  & . \\ . & . & . & . & . & {\cal
    R}_{2,1}^{2,1} & . & {\cal R}_{2,1}^{1,2} & .\\ . & . & {\cal
    R}_{0,2}^{2,0} & . & {\cal R}_{0,2}^{1,1} & . & {\cal
    R}_{0,2}^{0,2} & . & . \\ . & . & . & . & .  & {\cal
    R}_{1,2}^{2,1} & . & {\cal R}_{1,2}^{1,2} & .\\ . &. & . & .  &
  . & . & .  & . & {\cal R}_{2,2}^{2,2}
  \end{pmatrix}}
\end{equation}
{Indeed the `shape' - the non-vanishing pattern - is the same for $R$, $R_{rlex}$ and $\R$.}
The grlex matrix is obtained from this with
\[
R_{grlex}=P_G R_{rlex} P_G,
\]
where $P_G$ implements the transpositions $\ket{01}\leftrightarrow \ket{20}$ and $\ket{21}\leftrightarrow \ket{02}$.
Then an ACC matrix takes the block form exemplified by
\begin{equation}\label{eq:matf}
  {R_{grlex}=
\begin{pmatrix}{\cal R}_{0,0}^{0,0} &. & . & .
 & . & . & .  & . & .\\ . & {\cal R}_{1,0}^{1,0} & {\cal
    R}_{1,0}^{0,1} & . & . & . & . & . & .\\ . & {\cal R}_{0,1}^{1,0}
  & {\cal R}_{0,1}^{0,1} & . & .  & . & . & . & . \\ . & . & . & {\cal
    R}_{2,0}^{2,0} & {\cal R}_{2,0}^{1,1} &{\cal R}_{2,0}^{0,2 } & .
  & . & . \\ . & . & . & {\cal R}_{1,1}^{2,0} & {\cal R}_{1,1}^{1,1} &
 {\cal R}_{1,1}^{0,2} & . &  .  & . \\ . & . & . & {\cal
    R}_{0,2}^{2,0} & {\cal R}_{0,2}^{1,1} & {\cal
    R}_{0,2}^{0,2} & . & . & . \\ . & . & . & . & . & . & {\cal
    R}_{2,1}^{2,1} & {\cal R}_{2,1}^{1,2} & .\\ . & . & . & . & .  & . & {\cal
    R}_{1,2}^{2,1} & {\cal R}_{1,2}^{1,2} & .\\ . &. & . & .  &
  . & . & .  & . & {\cal R}_{2,2}^{2,2}
  \end{pmatrix}}
\end{equation}
In order to save space we will in the following
just give the blocks as
\beq\label{eq:Rblock}
R_{grlex}=\begin{bmatrix} {\cal R}_{0,0}^{0,0} \end{bmatrix}
  \begin{bmatrix}
  {\cal R}_{1,0}^{1,0} & {\cal R}_{1,0}^{0,1}\\ {\cal R}_{0,1}^{1,0} &
  {\cal R}_{0,1}^{0,1}\end{bmatrix}
\begin{bmatrix}   {\cal
    R}_{2,0}^{2,0} & {\cal R}_{2,0}^{1,1} &{\cal R}_{2,0}^{0,2 }
  \\ {\cal R}_{1,1}^{2,0} & {\cal R}_{1,1}^{1,1} & {\cal
    R}_{1,1}^{0,2} \\ {\cal R}_{0,2}^{2,0} & {\cal R}_{0,2}^{1,1} &
     {\cal R}_{0,2}^{0,2}\end{bmatrix}
  \begin{bmatrix}       {\cal
    R}_{2,1}^{2,1} & {\cal R}_{2,1}^{1,2} \\  {\cal
      R}_{1,2}^{2,1} & {\cal R}_{1,2}^{1,2} 
  \end{bmatrix}
  \begin{bmatrix} {\cal R}_{2,2}^{2,2} \end{bmatrix}.
\eeq
Recall that $\R$ is obtained from $R$ by exchanging lower indices,
which corresponds to up-down reflection within the block. In order to
match with \cite{MR} (using shifted basis labels $\{0,1,2\} \leadsto$
$\{1,2,3\}$),
highlight the new parameters, and 
save from writing
many
double indices we introduce
shorthand notation for $\R$:

\newcommand{\aaa}{{a}}
\newcommand{\bbb}{{b}}
\newcommand{\ccc}{{c}}
\newcommand{\ddd}{{d}}
\newcommand{\xone}{{x_1}}
\newcommand{\xfour}{{x_4}}

\begin{equation}\label{Rcheck}
\R =PR = 
\pmatrixx{
 a_{1}&\cdot&\cdot&\cdot&\cdot&\cdot&\cdot&\cdot&\cdot\cr 
 \cdot&\aaa_{12}&\cdot&\bbb_{12}&\cdot&\cdot&\cdot&\cdot&\cdot\cr 
 \cdot&\cdot&\aaa_{13}&\cdot&\xone&\cdot&\bbb_{13}&\cdot&\cdot\cr 
 \cdot&\ccc_{12}&\cdot&\ddd_{12}&\cdot&\cdot&\cdot&\cdot&\cdot\cr 
 \cdot&\cdot&x_{2}&\cdot&a_{2}&\cdot&x_{3}&\cdot&\cdot\cr 
 \cdot&\cdot&\cdot&\cdot&\cdot&\aaa_{23}&\cdot&\bbb_{23}&\cdot\cr 
 \cdot&\cdot&\ccc_{13}&\cdot&\xfour&\cdot&\ddd_{13}&\cdot&\cdot\cr 
 \cdot&\cdot&\cdot&\cdot&\cdot&\ccc_{23}&\cdot&\ddd_{23}&\cdot\cr 
 \cdot&\cdot&\cdot&\cdot&\cdot&\cdot&\cdot&\cdot&a_{3}\cr  
 }
 \end{equation}

  Then the block form is \beq\label{eq:blocsh}
\begin{bmatrix} a_{1}\end{bmatrix}
\begin{bmatrix}  a_{12} & b_{12} \\ c_{12} & d_{12} \end{bmatrix}
\begin{bmatrix} a_{13} & x_1 & b_{13} \\ x_{2} &  a_{2} & x_3\\
 c_{13} & x_4 & d_{13} \end{bmatrix}
\begin{bmatrix}  a_{23} & b_{23} \\ c_{23} & d_{23} \end{bmatrix}
\begin{bmatrix} a_{3} \end{bmatrix}
\eeq

\subsection{Symmetries} \label{ss:syms}

Naturally it is useful to consider additive charge-conserving solutions to \eqref{eq:braidr} up to  transformations that preserve \eqref{eq:braidr} and the additive charge-conserving condition. 

\begin{enumerate}
\item  \textbf{Scaling} symmetry: Equation  \eqref{eq:braidr}  and the additive charge-conserving condition is invariant under rescaling $\R$ by a non-zero complex number.
\item \textbf{Transpose} symmetry: The additive charge-conservation is preserved under
transpose: $\R\mapsto\R^T$;
and
of course
\eqref{eq:braidr} is satisfied by $\R^T$ if it is satisfied by $\R$
quite generally. Indeed, notice that from the form of \eqref{eq:YB} it is easy to see that if ${\cal
  \R}_{i,j}^{k,l}$ solves the equation, then so does ${\cal
  \R}_{k,l}^{i,j}$.
The effect on the variable choices in \eqref{Rcheck} are $b_{ij}\leftrightarrow c_{ij}$ and  $x_1\leftrightarrow x_2$ and $x_3\leftrightarrow x_4$.  
\item \textbf{Left-Right} (LR) symmetry: 
Changing the ordering of the basis from lex to rlex 
the resulting matrix will also be a solution. This can be seen in matrix entries because if ${\cal \R}_{i,j}^{k,l}$ solves
 equation \eqref{eq:braidr}, then so does ${\cal \R}_{j,i}^{l,k}$.  This corresponds to reflecting each of the blocks in \eqref{eq:blocsh} across both the diagonal and the skew-diagonal, i.e. $a_{ij}\leftrightarrow d_{ij}$, $b_{ij}\leftrightarrow c_{ij}$ as well as $x_1\leftrightarrow x_4$ and $x_2\leftrightarrow x_3$.  
\item \textbf{02-} or $\ket{0}\leftrightarrow\ket{2}$ symmetry:
while \eqref{eq:braidr} is clearly invariant under local basis changes, the additive charge-conserving condition is not.
However, the local basis change (permutation) $\ket{j}\leftrightarrow \ket{2-j}$ with indices $\{0,1,2\}$ taken modulo $3$ does preserve the form of an additive charge-conserving matrix: the span of the $\ket{ij}$ with $i+j=2$ is preserved, while the $\ket{ij}$ with $i+j=1$ and $i+j=3$ are interchanged as are the vectors $\ket{00}$ and
$\ket{22}$.
The effect on the block form \eqref{eq:blocsh} is to interchange the pairs of $1\times 1$ and $2\times 2$ blocks followed by a reflection across both the diagonal and skew-diagonal of each block.  \end{enumerate}

Of course these symmetries can be composed with one another and,
discounting the rescaling,
one
finds that the group of such symmetries is the dihedral group of order $8$.
This can be seen by tracking the orbit of the $2\times 2$ matrix $\begin{bmatrix}  a_{12} & b_{12} \\ c_{12} & d_{12} \end{bmatrix}$, since there are no symmetries that fix it.  Indeed, we see that there are $4$ forms it can take, generated by the reflections across the diagonal and, independently across the skew-diagonal, and two positions in \eqref{eq:blocsh} it can occupy.

\section{The solutions}  \label{S:sol}

\newcommand{\A}{A_R} 

For constant Yang--Baxter solutions, a necessary and sufficient set of 
constraint equations on the indeterminate matrix entries arise as follows.
Firstly compute, say,
\beq  \label{eq:AR}
\A := \R_1 \R_2 \R_1 - \R_2 \R_1 \R_2 
\eeq
which we call the \emph{braid anomaly} so that the constraints are obtained from $\A =0$.

The SCC case in which all $x_i$ vanish was solved in \cite{MR}.  Note that in ACC
some $x_i$ can be nonzero there will be mixing between more states $|ij\rangle$, but always with $i+j$ constant, so this is a
computationally relatively 
modest generalisation.
However the full symmetry of indices that
exists for SCC is now broken.
This ansatz-relaxing obviously
increases the complexity of the system of cubic equations, but they
can still be solved, as given below.

We organise the solutions according to which $x_i$s are
vanishing.  In principle there are $2^4-1=15$ cases (excluding the SCC
case), but we can use the above symmetries in order to omit some $x$
configurations.  This leads to the following classification into 6
cases:
\begin{enumerate}
\item All $x_i$ are nonzero. See \S\ref{ss:Jxxxx1} and \S\ref{ss:Jxxxx}.
\item Precisely one $x$ vanishes, by symmetry it can be assumed to be $x_4$.
  See
  \S\ref{ss:nosoln}.
\item $x_3x_4\neq0$ and $x_1=x_2=0$, related by the LR symmetry to
   $x_1x_2\neq0$ and $x_3=x_4=0$.
   See \S\ref{ss:case3}.
\item $x_1x_3\neq0$ and $x_2=x_4=0$, related to $x_2x_4\neq0$ and $x_1=x_3=0$ by
  transposition. See \S\ref{ss:nosoln}.
\item $x_1x_4\neq0$, $x_2=x_3=0$, related to $x_2x_3\neq0$ and $x_1=x_4=0$ by
  transposition, \S\ref{ss:case5}
\item Only one $x$ is nonzero, by symmetry it can be assumed to be $x_4$,
  \S\ref{ss:case6}.
  \end{enumerate}

As noted, solution of constant Yang--Baxter is equivalent to solving  $\A =0$.
We write out $\A$ explicitly in Appendix~\ref{ss:cubics}.
We solve for the various cases as above in the following
sections \ref{ss:Jxxxx1}--\ref{ss:case6}.
In the first of these we treat Case~1 relatively gently. After that we will
proceed more rapidly though all cases.

\subsection{The $\xxxxx_1 \xxxxx_2 \xxxxx_3 \xxxxx_4\neq 0$ solutions} \label{ss:Jxxxx1}

\newcommand{\ab}{a_2}

Recall the ACC ansatz for $\R$, which is as in  (\ref{Rcheck}).
Consider now the refinement of this ansatz indicated by the block structure

$$
\left[\begin {array}{c} 1\end {array} \right] \left[ 
\begin {array}{cc} 1&\cdot\\ \noalign{\medskip}\cdot&1\end {array} \right] 
 \left[ \begin {array}{ccc} a&x_{{1}}&b\\ \noalign{\medskip}{\frac {x_
{{3}} \left( a-1 \right) }{b}}&{\frac {x_{{1}}x_{{3}}+b}{b}}&x_{{3}}
\\ \noalign{\medskip}{\frac {{x_{{1}}}^{2}{x_{{3}}}^{2}}{{b}^{3}}}&-{
\frac {x_{{1}} \left( ab+x_{{1}}x_{{3}} \right) }{a{b}^{2}}}&-{\frac {
x_{{1}}x_{{3}}}{ab}}\end {array} \right] \left[ \begin {array}{cc} 1
&\cdot\\ \noalign{\medskip}\cdot&1\end {array} \right]  \left[ \begin {array}
     {c} 1\end {array} \right]
$$
that is

\beq  \label{eq:JaR1}
\R_{j} \; 
=
\left[
\begin{array}{ccccccccc}
1 & \cdot & \cdot & \cdot & \cdot & \cdot & \cdot & \cdot & \cdot 
\\
 \cdot & 1 & \cdot & \cdot & \cdot & \cdot & \cdot & \cdot & \cdot 
\\
 \cdot & \cdot & a  & \cdot & {\xxxxx_1}  & \cdot & b  & \cdot & \cdot 
\\
 \cdot & \cdot & \cdot & 1 & \cdot & \cdot & \cdot & \cdot & \cdot 
\\
 \cdot & \cdot & \frac{{\xxxxx_3} \left(a -1\right)}{b} & \cdot & \frac{{\xxxxx_1} {\xxxxx_3} +b}{b} & \cdot & {\xxxxx_3}  & \cdot & \cdot 
\\
 \cdot & \cdot & \cdot & \cdot & \cdot & 1 & \cdot & \cdot & \cdot 
\\
 \cdot & \cdot & \frac{{\xxxxx_3}^{2} {\xxxxx_1}^{2}}{b^{3}} & \cdot & -\frac{{\xxxxx_1} \left(a b +{\xxxxx_1} {\xxxxx_3} \right)}{a \,b^{2}} & \cdot & -\frac{{\xxxxx_3} {\xxxxx_1}}{a b} & \cdot & \cdot 
\\
 \cdot & \cdot & \cdot & \cdot & \cdot & \cdot & \cdot & 1 & \cdot 
\\
 \cdot & \cdot & \cdot & \cdot & \cdot & \cdot & \cdot & \cdot & 1 
\end{array}\right]
\eeq 
Here the parameters 
$\aaaaa_{13}, \;  \bbbbb_{13}, \; \xxxxx_1, \; \xxxxx_3$ 
are indeterminate 
(we write $\aaaaaa = \aaaaa_{13}, \; \bbbbbb = \bbbbb_{13}$)
but  the remaining parameters are replaced
with functions of these four as shown.

\begin{proposition}   (I) Consider the ansatz for $\R$ in (\ref{Rcheck}).
If we leave 
parameters
$\aaaaa_{13}, \; 
\bbbbb_{13}, \; \xxxxx_1, \; \xxxxx_3$ 
indeterminate 
(here we write $\aaaaaa = \aaaaa_{13}, \;
 \bbbbbb = \bbbbb_{13}$)
but replace the 
remaining parameters 
with functions of these four as
shown in (\ref{eq:JaR1})
then the braid anomaly $\A$ has an overall factor
\[
\x_3^{2} \x_1^{2} a 
+ a^{2} b^{2}
+ \x_1  b  \x_3 a 
- a \,b^{2}
- \x_1 b \x_3
\;\; = \;\;
b^2 a (a-1) + b \x_3 \x_1 (a-1) + a \x_3^2 \x_1^2
\]
That is, we have a family of solutions 
obeying $\R_1 \R_2 \R_1 = \R_2 \R_1 \R_2$
with 
free 
non-zero
parameters (say) $a, \x_1, \x_3$, and parameter
$b$ determined by 
\beq \label{eq:boxx} 
\frac{b}{\x_1 \x_3} = 
\frac{-\frac{1}{a} \pm \sqrt{\frac{1}{a^2}-\frac{4}{a-1}}}{2}
\hspace{1cm}\mbox{ or } \hspace{1cm}
\frac{\x_1 \x_3}{b} = 
\frac{-(a-1)\pm\sqrt{(a-1)^2 -4a^2 (a-1)}}{2a}
\eeq 
and the remaining entries determined as 
in (\ref{eq:JaR1}) above. 

\noindent (II) If $\x_1 \x_2 \x_3 \x_4 \neq 0$ then the above  
(with $a\neq 1$)
gives the complete set of solutions
up to overall rescaling. 
\end{proposition}

\newcommand{\SSSSS}{{\mathsf{S}}_{R}}

\noindent {\em Proof.}
(I) is simply a brutal but straightforward calculation, plugging in to
$\A$ as given for example in Appendix \ref{ss:cubics}.
$\;$ 
For (II) we proceed as follows. 
The matrix $\A$ is rather large to write out
 (again see Appendix \ref{ss:cubics}),
  but a subset of its entries is
\beq \label{eq:SR0} 
\SSSSS \;\; =\;\;   \{ \;
-\aaaaa_{12} \bbbbb_{12} \ccccc_{12}
-\aaaaa_1 \aaaaa_{12} (\aaaaa_{12} - \aaaaa_{1}),
\;\;\;\;\;\;\; \hspace{.1cm} 
- \bbbbb_{23} \ccccc_{23} \ddddd_{23}
-\aaaaa_3 \ddddd_{23} (\ddddd_{23} - \aaaaa_{3}),
\;\;\;\;\;
\hspace{1cm} 
\eeq 
\beq   \label{eq:SR00}   \hspace{1.2cm}
-\bbbbb_{23} \xxxxx_1 \xxxxx_3,
\;\;\;\;\;\;\;\hspace{.1cm}
-\bbbbb_{12} \xxxxx_1 \xxxxx_3,
\;\;\;\;\;\;
(\aaaaa_{12}-\ddddd_{12}) \xxxxx_1 \xxxxx_3,
\;\;\;\;\;\;
((-\ddddd_{12} +\aaaaa_1 ) \bbbbb_{13} -\bbbbb_{12}^2 ) \xxxxx_1,
\;\;\;\;\;
\eeq 
\[
\aaaaa_{13} (\bbbbb_{12} \ccccc_{12} 
-\bbbbb_{23} \ccccc_{23})
+\aaaaa_{12} \aaaaa_{23}
(\aaaaa_{12} - \aaaaa_{23}),
\hspace{1cm}
(\ddddd_{12} -\ddddd_{23} ) \xxxxx_1 \xxxxx_2,
\hspace{1cm}
(\ddddd_{23} -\aaaaa_{23} ) \xxxxx_1 \xxxxx_3,
\]
\beq  \label{eq:SR}
\hspace{-.1cm} 
-\aaaaa_{12} \xxxxx_1 \xxxxx_3 -\aaaaa_{13} \bbbbb_{13} \ddddd_{13},
\;\;\;\;\;\hspace{1cm} 
-\aaaaa_{12} \xxxxx_2 \xxxxx_4
-\aaaaa_{13} \ccccc_{13} \ddddd_{13},
\hspace{.9cm}
\aaaaa_{12} \ccccc_{12} \ddddd_{12},
\;\;\;\;\;\hspace{1cm} 
\aaaaa_{23} \ccccc_{23} \ddddd_{23},
\;\;\;\;\;\; 
\eeq
\beq  \label{eq:SR4}
-\aaaaa_{13} \bbbbb_{13} \xxxxx_4
+(\aaaaa_1 \aaaaa_{12}
-\aaaaa_{12} \aaaaa_2
-\aaaaa_1 \aaaaa_{13} ) \xxxxx_1,
\eeq
\beq  \label{eq:SR5}
( \aaaaa_1 \ddddd_{13}
+\aaaaa_2 \ddddd_{12}
-\aaaaa_1 \ddddd_{12} ) \xxxxx_3
+\bbbbb_{13} \ddddd_{13} \xxxxx_2,
\hspace{1.2cm} 
(\aaaaa_{13} \ddddd_{13}
+\aaaaa_2 \aaaaa_{23}
-\aaaaa_{13} \aaaaa_{23} ) \xxxxx_3
+ ( \aaaaa_3 \bbbbb_{13} 
-\bbbbb_{23}^2 ) \xxxxx_2,
\;\;\;\;\;\;  
\eeq
\beq \label{eq:SR7}
\hspace{.03cm} 
\aaaaa_{13} \ccccc_{13} \xxxxx_3
+(\aaaaa_{13} \aaaaa_3
-\aaaaa_{23} \aaaaa_3
+\aaaaa_2 \aaaaa_{23} ) \xxxxx_2,
\hspace{1.1cm}
\aaaaa_{12} \ccccc_{13} \xxxxx_3
+(-\aaaaa_2 \ddddd_{23}
+\aaaaa_{13} \ddddd_{23}
-\bbbbb_{23} \ccccc_{12}
+\aaaaa_2^2 ) \xxxxx_2,
\eeq
\beq \label{eq:SR77}
( \aaaaa_1 \ccccc_{13} -\ccccc_{12}^2 ) \xxxxx_3
+ ( \aaaaa_{13} \ddddd_{13}
-\ddddd_{12} \ddddd_{13}
+\aaaaa_2 \ddddd_{12} ) \xxxxx_2,
\eeq 
\beq \label{eq:SR777}  
-\aaaaa_2 \xxxxx_1 \xxxxx_2
+\aaaaa_{13} \ddddd_{12}^2
-\aaaaa_{13}^2 \ddddd_{12}
-\aaaaa_{23} \bbbbb_{13} \ccccc_{13}
+\aaaaa_{23} \bbbbb_{12} \ccccc_{12},
\eeq
\beq \label{eq:SR8}
\aaaaa_1 \xxxxx_3 \xxxxx_4
+\ddddd_{13} \xxxxx_1 \xxxxx_2
-\aaaaa_2 \ddddd_{12}^2
-\bbbbb_{12} \ccccc_{12} \ddddd_{12}
+\aaaaa_2^2 \ddddd_{12},
\;\;\;\;\hspace{.42cm} 
\eeq
\beq \label{eq:SR9}
\aaaaa_{13} \xxxxx_3 \xxxxx_4
+\aaaaa_3 \xxxxx_1 \xxxxx_2
-\aaaaa_{23} \bbbbb_{23} \ccccc_{23}
-\aaaaa_2 \aaaaa_{23}^2
+\aaaaa_2^2 \aaaaa_{23},
\;\;\;\;
...
\} \;
\eeq 
Imposing $\A=0$ 
and  $\xxxxx_1 \xxxxx_3 \neq 0$ 
we thus get $\bbbbb_{12} = \bbbbb_{23}=0$
and $\aaaaa_{12} = \ddddd_{12}$ from (\ref{eq:SR00}).
Note that $\R$ is invertible, so 
$\aaaaa_1 , \aaaaa_3  \neq 0$
and 
\beq \label{eq:Rinv}
\aaaaa_{12} \ddddd_{12}-\bbbbb_{12} \ccccc_{12} \neq 0 ,
\hspace{1cm} 
\aaaaa_{23} \ddddd_{23}-\bbbbb_{23} \ccccc_{23} \neq 0.
\eeq 
Thus 
$
\aaaaa_{12} , \ddddd_{12},\aaaaa_{23}, \ddddd_{23} \neq 0
$. 
Thus $\ddddd_{12}=\aaaaa_{23} =\aaaaa_{12} = \aaaaa_1
= \ddddd_{23} = \aaaaa_3$
and $\ccccc_{12} = \ccccc_{23}=0$
from (\ref{eq:SR}).
Note that if $\R$ is a solution then so is any non-zero scalar multiple, so 
we first scale $\R$ by an overall factor, so that 
$\aaaaa_1 =1$.
This confirms the form for $\R$ above outside the 
3x3 block. 

Observe now from (\ref{eq:SR}) that 
$ \aaaaa_{13} \bbbbb_{13} \ddddd_{13} 
=-\xxxxx_1 \xxxxx_3 \neq 0$,
and that we may either replace 
$ \aaaaa_{13} = 
-\frac{\xxxxx_1 \xxxxx_3}{\bbbbb_{13} \ddddd_{13}} $,
or 
$ \ddddd_{13} = 
-\frac{\xxxxx_1 \xxxxx_3}{\aaaaa_{13} \bbbbb_{13}}$.
The latter gives the form of $\ddddd_{13}$ in the Proposition.

Before proceeding we will need to show that 
$\dddddd=1$ cannot occur here.
(Recall $\aaaaaa=\aaaaa_{13}$, and write also $\dddddd=\ddddd_{13}$.)
Comparing (\ref{eq:SR8}) and (\ref{eq:SR9}) we find
\[
(\aaaaaa-1) \xxxxx_3 \xxxxx_4 
= (\dddddd-1) \xxxxx_1 \xxxxx_2 .
\]
So $\aaaaaa=1$ if and only if $\dddddd=1$. 
So if $\aaaaaa=1$ then 
$\bbbbb_{13} =-\xxxxx_1 \xxxxx_3$ and 
$\ccccc_{13} = -\xxxxx_2 \xxxxx_4$
(consider (\ref{eq:SR}i/ii)).

Evaluating (\ref{eq:SR7}ii)-(\ref{eq:SR77}) here we find:
\[
\aaaaa_2^2 -\aaaaa_2 +a -\aaaaa_2 +\dddddd -a \dddddd \;=\;\; 
(\aaaaa_2 -1)^2 -(a-1)(\dddddd-1) \; = 0
\]
So if $\aaaaaa=1$ then $\aaaaa_2=1$.
Further, if $\aaaaaa=1$ then (\ref{eq:SR777}) becomes
$
-\xxxxx_1 \xxxxx_2 -\bbbbb_{13} \ccccc_{13} = -\xxxxx_1 \xxxxx_2 - \xxxxx_1 \xxxxx_2 \xxxxx_3 \xxxxx_4 =0
$
so $\xxxxx_3 \xxxxx_4 = -1$. 
But if $\aaaaaa=1$ then (18) gives 
$\xxxxx_4 = (-\xxxxx_1)/(-\xxxxx_1 \xxxxx_3) = 1/\xxxxx_3$ -- a contradiction. 
We conclude here that $(\aaaaaa-1)\neq 0$;  and hence 
$(\dddddd-1) \neq 0$.

From (\ref{eq:SR5}) we have two formulae for $\bbbbb_{13} \xxxxx_2$. Equating we have
\[
\frac{( \aaaaa_1 \ddddd_{13}
+\aaaaa_2 \ddddd_{12}
-\aaaaa_1 \ddddd_{12} )}  {\ddddd_{13} }  \hspace{.2cm} 
=
(\aaaaa_{13} \ddddd_{13}
+\aaaaa_2 \aaaaa_{23}
-\aaaaa_{13} \aaaaa_{23} ) 
\]
that is 
\[
\frac{(  \ddddd_{13}   +\aaaaa_2  - 1 )}  
{\ddddd_{13} }  
\hspace{.2cm} 
=
(\aaaaa_{13} \ddddd_{13}
+\aaaaa_2 
-\aaaaa_{13}  ) ,\;\;\text{thus}\;\;
\frac{(  \aaaaa_2  - 1 )(1- \ddddd_{13})}  
{\ddddd_{13} }  
\hspace{.2cm} 
=
-\aaaaa_{13} (1-\ddddd_{13} ).
\]
Since
$\ddddd_{13} -1 \neq 0$ 
we have
$
 \aaaaa_2  - 1   =  -\aaaaa_{13} \ddddd_{13}
 = \frac{ \xxxxx_1 \xxxxx_3 }{ \bbbbb_{13} } 
$
giving $ \aaaaa_2$ as in the Proposition. 
Plugging back in 
we find
$(\ddddd_{13} +\aaaaa_2 -1) \xxxxx_3 
=-\bbbbb_{13} \ddddd_{13} \xxxxx_2,
\;\;
\frac{(\aaaaa_{13} -1) \xxxxx_3}{\bbbbb_{13}} = \xxxxx_2
$
as in the Proposition.

From (\ref{eq:SR4}) we have
\[
\xxxxx_4 \; = \; 
\xxxxx_1 \frac{1-\aaaaa_2 -\aaaaa_{13}}{\aaaaa_{13} \bbbbb_{13}}
\; = \; 
\xxxxx_1 \frac{-\xxxxx_1 \xxxxx_3 -\aaaaa_{13} \bbbbb_{13}}{\aaaaa_{13} \bbbbb_{13}^2} 
\]
as in the Proposition. 
Finally from (\ref{eq:SR7}) we  now have
\[
\ccccc_{13} = 
-\frac{ ( \aaaaa_{13} +\aaaaa_2 -1 ) \xxxxx_2}{\aaaaa_{13} \xxxxx_3}
=
-\frac{(\aaaaa_{13} \bbbbb_{13} 
+\xxxxx_1 \xxxxx_3 )(\aaaaa_{13}-1)}{\aaaaa_{13} \bbbbb_{13}^2}
\]
\[
 \ccccc_{13} \xxxxx_3
= - ( \aaaaa_{13} 
+\aaaaa_2 (\aaaaa_2 -1) ) \xxxxx_2
= - ( {a}+ \frac{( \xxxxx_1 \xxxxx_3 +b )}{b}  
\frac{\xxxxx_1 \xxxxx_3}{ b}  )  
\frac{ (\aaaaa_{13}-1) \xxxxx_3 }{b}\]
\[= \frac{ {-}
\xxxxx_3 (
(a b^2 + \xxxxx_1 \xxxxx_3 (\xxxxx_1 \xxxxx_3 +\bbbbb_{13}))(\aaaaa_{13}-1))}{\bbbbb_{13}^3}
\]

Equating the two formulae for $\ccccc_{13}$,
and noting that $\aaaaa_{13}-1\neq 0$, we have
\[
\aaaaa_{13} (\xxxxx_1 \xxxxx_3 )^2 
+(\aaaaa_{13}-1) \bbbbb_{13} \xxxxx_1 \xxxxx_3
+(\aaaaa_{13}-1) \aaaaa_{13} \bbbbb_{13}^2 
=0
\]
Plugging back in to (\ref{eq:SR7}) 
we obtain $\ccccc_{13}$ as in the 
Proposition, so we are done.
\hfill 
\qed

\newcommand{\mycomment}[1]{}
\newcommand{\JHm}[1]{\textcolor{green}{#1}}
\allowdisplaybreaks

\subsection{Case 1: $x_1x_2x_3x_4\neq0$ revisited \label{ss:Jxxxx}}

In this section we solve the case in which $x_1x_2x_3x_4\neq0$
 again, but leaning directly on the Appendix
 (as we shall below for the remaining cases).
Since all $x_i$ are nonzero, we
conclude from equations \eqref{A5}-\eqref{A8} that
$b_{12}=c_{12}=b_{23}=c_{23}=0$, and from \eqref{A11}-\eqref{A13}
$a_{12}=d_{12}=a_{23}=d_{23}$. Then since $a_{12}\neq0$ from
\eqref{A29} we get $a_{12}=1$ and since $a_3\neq0$ from \eqref{A30}
$a_3=1$.

Now from \eqref{A18} we find $a_{13}b_{13}d_{13}\neq0$ and we can solve
$
d_{13}=-\frac{x_1x_3}{a_{13}b_{13}},
$
and from \eqref{A23}
$
c_{13}=\frac{b_{13}x_2x_4}{x_1x_3}.
$
Then from \eqref{A70} we find
$
x_4=\frac{x_1(1-a_2-a_{13})}{a_{13}b_{13}}.
$
Now it turns out that some equations factorize, for example
\eqref{A94} can be written as
$
x_1(a_{13}-1)[(a_2-1)b_{13}-x_1x_3]=0.
$
If we were to choose $a_{13}=1$ we reach a contradiction: from
\eqref{A58} we get $a_2=1$ and then \eqref{A50} and \eqref{A66} are
contradictory since $x_i\neq0$. Thus we can solve
$
a_2=\frac{b_{13}+x_1x_3}{b_{13}},
$
and then from \eqref{A50}
$
x_2=\frac{x_3(a_{13}-1)}{b_{13}}.
$

After this all remaining nonzero equations simplify to
\[
b_{13}^2a_{13}(a_{13}-1)+b_{13} x_1x_3(a_{13}-1)+a_{13}x_1^2x_3^2=0.
\]
This biquadratic equation can be resolved using Weierstrass elliptic
function $\wp$:\\
\[ a=-\wp+\frac5{12},\quad \beta=6\frac{12\wp+7+12\wp'}{(12\wp-5)(12\wp+7)},\quad
(\wp')^2=4\wp^3- \frac1{12}\wp+\frac{ 7\cdot  23}{2^3\ 3^3}
\]
where
$a_{13}=a$ and $b_{13}=x_1x_3\b$. The solution in
block form now reads
\bse\label{sol:case1}
\beq [1]  \begin{bmatrix}1&.\cr .&1\end{bmatrix}  
\begin{bmatrix}a&x_{1}&\b 
 x_{1}  x_{3}\cr \frac{a -1}{ \b  x_{1}}& \frac{\b +1}{
 \b}&x_{3}\cr \frac{1}{ \b^{3}  x_{1}  x_{3}}& \frac{- (a 
 \b +1 ) }{ a  \b^{2}  x_{3}}& \frac{-1}{ a  \b}\end{bmatrix} 
 \begin{bmatrix}1&.\cr .&1\end{bmatrix}  [1]
\eeq
with constraint
\beq
\b^2 a (a - 1) + \b (a - 1) + a=0.
\eeq
\ese

\subsection{Cases 2 and 4:  $x_1x_3\neq0$ and $x_4x_2=0$ \label{ss:nosoln}}

From \eqref{A7},\eqref{A8} we get $b_{12}=b_{23}=0$ and then since the
matrix is non-singular we must have $a_{12}d_{12}a_{23}d_{23}\neq 0$.
Then from \eqref{A1},\eqref{A2} we get $c_{12}=c_{23}=0$, from
\eqref{A29},\eqref{A31} $a_{12}=d_{12}=1$ (recall that we have scaled $a_1=1$) and from
\eqref{A84},\eqref{A85} $a_{23}=d_{23}=1$. Next from \eqref{A30}
$a_3=1$.  Since $x_1x_3\neq0$ we have from \eqref{A18} that
$a_{13}b_{13}d_{13}\neq 0$ and then from \eqref{A24} we find
$c_{13}=0$.

To continue we consider first the case $x_2=0$, $x_4$ free. Then from
\eqref{A52} we get $a_{13}=1$ and from \eqref{A76} $d_{13}=1-a_2$.
Then since $d_{13}\neq0$ we cannot have $a_2=1$ but \eqref{A106} is
$x_3(a_2-1)^2=0$, a contradiction.

Next assume $x_2\neq0,\,x_4=0$. Then from \eqref{A66} we get
$d_{13}=1$ and from \eqref{A70} $a_{13}=1-a_2$ but \eqref{A98}
yields $a_2=1$ which is in contradiction with $a_{13}\neq0$.

Thus there are no solutions in this case.

\subsection{Case 3: $x_1=x_2=0$, $x_3x_4\neq0$ \label{ss:case3}}

From \eqref{A16} we get $a_{23}=a_{12}$ and from \eqref{A94} and
\eqref{A104} $b_{13}=b_{12}^2$ and $c_{13}=c_{12}^2$.  On the basis of
\eqref{A42} and \eqref{A46} we can divide the problem into 2 branches:
Case 3.1: $a_{12}=0$, $b_{12}c_{12}\neq 0$, and
Case 3.2: $a_{12}d_{12}\neq0$, $b_{12}=c_{12}= 0$.

\paragraph{Case 3.1: $a_{12}=0$, $b_{12}c_{12}\neq 0$.}

From \eqref{A46} we get $a_{13}=0$ and then from \eqref{A41}
$d_{12}=1-b_{12}c_{12}$ and from \eqref{A76}
$d_{13}=(1-a_2)(1-b_{12}c_{12})$. Then from \eqref{A100} and
\eqref{A101}, $b_{23}=a_2^2/c_{12}$, $c_{23}=a_2^2/b_{12}$ and from
\eqref{A95} $a_3=a_2^4/(b_{12}c_{12})^2$. Since $a_3\neq0$ we have
$a_2\neq0$ and can solve $d_{23}$ from \eqref{A43}:
$d_{23}=(a_2^4-(b_{12}c_{12})^3)/(b_{12}c_{12})^2$.

Now \eqref{A82} factorizes as
$(a_2^2-b_{12}c_{12})(a_2^2+a_2b_{12}c_{12}+(b_{12}c_{12})^2)=0$.

Case 3.1.1: If we choose the first factor and set
$b_{12}=a_2^2/c_{12}$ the remaining equations simplify to
$x_3=a_2(a_2^2-1)^2/x_4$, yielding the first solution ($a_2\to
a,\,c_{12}\to c$);
\beq \label{eq:III.1.1}
    [1] \begin{bmatrix}.& \frac{a^{2}}{ c}\cr c&
  1-a^2 \end{bmatrix} \begin{bmatrix}.&.& \phantom{c_{I}} \frac{a^{4}}{ c^{2}}
  \\
  .&a& \frac{a(a^2-1)^2}{x_4}  \\ c^{2}&x_{4}& (a +1 ) (a -1 )
  ^{2}\end{bmatrix} \begin{bmatrix}.& \frac{a^{2}}{ c}\cr c&
  1-a^2 \end{bmatrix} [1]
\eeq

\noindent {The eigenvalues of this solution are $1,-a^2$ and $a^3$ with multiplicities $5,3$ and $1$, respectively.}
\medskip

Case 3.1.2: For the second solution we solve \eqref{A82} by
$b_{12}=a\omega/c_{12}$, where $\omega$ is a cubic root of unity
$\omega\neq1$. Then the remaining equation is solved by
$x_3=a_2(a_2-1)(1-\omega a_2)/x_4$ and we have
\beq \label{eq:III.1.2}
    [1]  \begin{bmatrix}.& \frac{\omega   a}{ c}\\ c&
1-\omega a \end{bmatrix}
\begin{bmatrix}.&.& \frac{ \omega^2 a^{2} }{ c^{2}}\\
  .&a& \frac{(\omega^2-a )   (a -1 )   a}{ x_{4}}\\
  c^{2}&x_{4}& (1-\omega a)(1-a) \end{bmatrix} 
\begin{bmatrix}.& \frac{a^{2}}{ c}\\ \omega^2  a  c  & \omega a (a-1)\end{bmatrix} [ \omega a^{2} ]
\eeq
\noindent {The eigenvalues are $\{1,-\omega a,\omega a^2\}$ each with multiplicity $3$.}

For both solutions $a\neq 0,1$ and for the first $a\neq -1$.  Note that for the second case if $a=-1$ there are only two eigenvalues: $1$ and $\omega$.

\paragraph{Case 3.2:} $a_{12}d_{12}\neq0$, $b_{12}=c_{12}= 0$. From equations \eqref{A9}, \eqref{A29} we get
$a_{12}=d_{12}=1$ and from \eqref{A44}, \eqref{A47}
$b_{23}=c_{23}=0$. Due to non-singularity we may now assume
$a_{23}d_{23}\neq0$ and then from  \eqref{A10}, \eqref{A16} we get
$a_{23}=d_{23}=1$. Since $a_3\neq0$ \eqref{A30} yields $a_3=1$.
Now from \eqref{A69} and \eqref{A100} we get $a_2=1,d_{13}=0$,
after which we get a contradiction in \eqref{A50}.

\subsection{Case 5:  $x_2x_3\neq0$, $x_1=x_4=0$ \label{ss:case5}}
This case contains many solutions and therefore it is necessary to do
some basic classification first. We do this on the basis of the
$2\times2$ blocks.

For the first $2\times2$ block (the ``12'' block), consider equations
\eqref{A1}, \eqref{A2}, \eqref{A9}, \eqref{A29} and \eqref{A31}. The
solutions to these equations can be divided into the following:
\begin{itemize}
\item[{$\alpha$:}] $a_{12}d_{12}\neq0$. Then one finds $b_{12}=c_{12}=0$ and
  $a_{12}=d_{12}=a_1$.
\item[{$\beta$:}] $a_{12}\neq0$, $d_{12}=0$ and $b_{12}c_{12}\neq0$, then
  $a_{12}=a_1-b_{12}c_{12}/a_1$
\item[{$\gamma$:}] $d_{12}\neq0$, $a_{12}=0$ and $b_{12}c_{12}\neq0$, then
  $d_{12}=a_1-b_{12}c_{12}/a_1$
\item[{$\delta$:}] $a_{12}=d_{12}=0$.
\end{itemize}
 The results for the other $2\times2$ block (the ``23'' block) are
 obtained by index changes, including $a_1 \to a_3$, we denote them as
 $\alpha'$ etc.

In principle there would be $4\times4=16$ cases, but we can omit
several using the known symmetries. First of all for the ``12'' block
we can omit $\gamma$ because it is related to $\beta$ by LR
symmetry. The list of cases is as follows:
\begin{enumerate}
\item $(\alpha,\alpha')$:
  $[a_1]\begin{bmatrix}a_1&.\\.&a_1\end{bmatrix}[3\times3]
  \begin{bmatrix}a_3&.\\.&a_3\end{bmatrix}[a_3]$.
  \vskip.1cm
\item $(\alpha,\beta')$:
  $[a_1]\begin{bmatrix}a_1&.\\.&a_1\end{bmatrix}[3\times3]
  \begin{bmatrix}a_3-b_{23}c_{23}/a_3&b_{23}\\c_{23}&.\end{bmatrix}[a_3]$.  \vskip.1cm
\item $(\alpha,\delta')$:
  $[a_1]\begin{bmatrix}a_1&.\\.&a_1\end{bmatrix}[3\times3]
  \begin{bmatrix}.&b_{23}\\c_{23}&.\end{bmatrix}[a_3]$.  \vskip.1cm
\item $(\beta,\beta')$:
  $[a_1]
    \begin{bmatrix}a_1-b_{12}c_{12}/a_1&b_{12}\\c_{12}&.\end{bmatrix}
  [3\times3]
  \begin{bmatrix}a_3-b_{23}c_{23}/a_3&b_{23}\\c_{23}&.\end{bmatrix}[a_3]$.  \vskip.1cm
\item $(\beta,\gamma')$:
   $[a_1]
    \begin{bmatrix}a_1-b_{12}c_{12}/a_1&b_{12}\\c_{12}&.\end{bmatrix}
  [3\times3]
  \begin{bmatrix}.&b_{23}\\c_{23}&a_3-b_{23}c_{23}/a_3\end{bmatrix}[a_3]$.  \vskip.1cm
\item $(\beta,\delta')$: $[a_1]
    \begin{bmatrix}a_1-b_{12}c_{12}/a_1&b_{12}\\c_{12}&.\end{bmatrix}
  [3\times3]
  \begin{bmatrix}.&b_{23}\\c_{23}&.\end{bmatrix}[a_3]$.  \vskip.1cm
\item $(\delta,\delta')$: $[a_1]
    \begin{bmatrix}.&b_{12}\\c_{12}&.\end{bmatrix}
  [3\times3]
  \begin{bmatrix}.&b_{23}\\c_{23}&.\end{bmatrix}[a_3]$. 
\end{enumerate}
Here we have omitted $(\alpha,\gamma')$, $(\beta,\alpha')$,
$(\delta,\alpha')$, $(\delta,\beta')$ and $(\delta,\gamma')$, because
they are related to entries in the above list of seven by some symmetry.
 {Specifically, notice that the vanishing of $x_1$ and $x_4$ and the non-vanishing of $x_2x_3$ is preserved under LR-symmetry and the $\ket{0}\leftrightarrow\ket{2}$ symmetry, but \emph{not} the transpose symmetry. Moreover the composition of the LR and 02-symmetries has the effect of simply interchanging the pairs of $2\times 2$ and $1\times 1$ blocks.}

  \paragraph{Case 5.1} $(\alpha,\alpha')$
  
We scale to $a_1=1$ and from \eqref{A30} get $a_3=1$. According to
\eqref{A86} $a_2^2-a_2=0$ and then from \eqref{A106} and \eqref{A108} we
get $d_{13}=-b_{13}x_2/x_3$ and $c_{13}=-c_{13}x_2/x_3$ but then the
$3\times3$ block matrix becomes singular. Therefore no solutions for
this subcase.

\paragraph{Case 5.2}  $(\alpha,\beta')$

From  \eqref{A39} and \eqref{A43} we get
$c_{13}=c_{23}^2/a_3$ and $b_{13}=b_{23}a_3/c_{23}$. Next from
\eqref{A58} $d_{13}=0$ and from \eqref{A76} $a_2=1$ and from
\eqref{A54} $a_{23}=a_{13}$. After setting $a_3=-x_3c_{23}^2/x_2$ from
\eqref{A104}, the GCD of the remaining equations is
$(x_3c_{23}-x_2b_{23})(x_2+x_3c_{23})^2$ and we get two solutions:
( $c_{23}\to c,\, b_{23}\to b$)

5.2.1: $x_2=-x_3c_{23}^2$
\beq \label{eq:V.2.1}
    [1]  \begin{bmatrix}1&.\\ .&1\end{bmatrix} \begin{bmatrix}1-b
  c&.& \frac{b}{ c}\\ -x_{3}  c^{2}&1&x_{3}\\
  c^{2}&.&.\end{bmatrix} \begin{bmatrix}1-b  c &b\\
    c&.\end{bmatrix} [1]
  \eeq
noindent{The eigenvalues are $-bc$ with multiplicity $2$ and $1$ with multiplicity $7$. }

5.2.2: $x_2=x_3c_{23}/b_{23}$
\beq\label{5.2.2} [1]  \begin{bmatrix}1&.\\ .&1\end{bmatrix} \begin{bmatrix}1-b  c
 &.& -b^{2}\\ \frac{x_{3}  c}{ b}&1&x_{3}\\
 \frac{-c}{ b}&.&.\end{bmatrix} \begin{bmatrix}1-b  c
 &b\\ c&.\end{bmatrix} [-b  c] \eeq

\medskip\noindent
{The eigenvalues are $-bc$ with multiplicity $3$ and $1$ with multiplicity $6$}.
\medskip

\paragraph{Case 5.3} $(\alpha,\delta')$

From \eqref{A54} and \eqref{A62} we get $a_{13}=d_{13}=0$ and from
\eqref{A72} $a_2=1$.  Next \eqref{A40} and \eqref{A43} yield
$c_{13}=b_{13}=a_3$ and \eqref{A102} $x_3=-x_2a_3$ The remaining
equations are satisfied with $a_3=\epsilon_1$ and $c_{23}=\epsilon_2$,
where $\epsilon_j^2=1$. The result is
\begin{equation*}
    [1] \begin{bmatrix}1&.\\ .&1\end{bmatrix} \begin{bmatrix}.&.&
        \epsilon_{1}\\ x_{2}&1& -x_{2}
        \epsilon_{1}\\ \epsilon_{1}&.&.\end{bmatrix} \begin{bmatrix}.&
        \epsilon_{2}\\ \epsilon_{2}&.\end{bmatrix}
      [\epsilon_{1}]
      \end{equation*}
\medskip\noindent
{The eigenvalues are $1$ and $-1$ with multiplicity $7$ and $2$ if $\epsilon_1=1$ and $6$ and $3$ otherwise.  However, when $\epsilon_1=-1$ this is a special case of \eqref{5.2.2} by setting $b=c=-1$.  For $\epsilon_1=1$ we may take $b=c=1$ in \eqref{eq:V.2.1}.  Thus this case may be discarded \emph{a posteriori} as a subcase. }
\medskip

      \paragraph{Case 5.4} $(\beta,\beta')$
      
From \eqref{A37} and \eqref{A39} we get $c_{13}=c_{12}^2$ and
$a_3=c_{23}^2/c_{12}^2$. Next since $a_{12}=1-b_{12}c_{12}\neq0$ we
get $d_{13}=0$ from \eqref{A61}. From \eqref{A41}
$b_{13}=b_{12}/c_{12}$ and from \eqref{A78}
$a_{13}=1-b_{12}c_{12}$. For nonsingularity we must have $a_2\neq0$
and then from \eqref{A82} we get $b_{23}=b_{12}$ and from \eqref{A43}
 $c_{23}=c_{12}$. Now from \eqref{A74} we find $a_2=-x_3
c_{12}^2/x_2$ and after that the remaining equations factorize and we
have two solutions:

5.4.a $x_2=-c_{12}^2x_3$
\beq \label{5.4.a}
    [1]  \begin{bmatrix}1-bc&b\\ c&.\end{bmatrix}
  \begin{bmatrix}1-bc&.& \frac{b}{
 c}\\ -x_{3} 
 c^{2}&1&x_{3}\\ c^{2}&.&.\end{bmatrix}
  \begin{bmatrix}1-bc&b\\ c&.\end{bmatrix} [1]
  \eeq
noindent
{Eigenvalues are $1$ with multiplicity $6$ and $-bc$ with multiplicity $3$.}

5.4.b $x_2=c_{12}x_3/b_{12}$

\beq \label{5.4.b}
    [1]  \begin{bmatrix}1-bc&b\\ c&.\end{bmatrix}
  \begin{bmatrix}1-bc&.& \frac{b}{ c}\\
 \frac{x_{3}  c}{ b}& -b  c&x_{3}\\
 c^{2}&.&.\end{bmatrix} \begin{bmatrix}1-bc&b\\
    c&.\end{bmatrix} [1]
  \eeq
\noindent
  {Eigenvalues are $1$ with multiplicity $5$ and $-bc$ with multiplicity $4$.}

\paragraph{Case 5.5} $(\beta,\gamma')$

Since the matrix is non-singular we must have $a_2\neq0$.  From
\eqref{A38} and \eqref{A41} we get $c_{13}=c_{12}^2$ and
$b_{13}=b_{12}/c_{12}$, and from \eqref{A39} and \eqref{A43}
$b_{12}=b_{23}^2c_{12}/a_3$, $c_{23}=b_{23}c_{12}^2/a_3$. Then we get
from several equations the condition $a_{13}d_{13}=0$. If both
$a_{13}=d_{13}=0$, we would get from \eqref{A78} $a_{12}=0$,
which would lead to case $\delta'$. Therefore we have two branches:

5.5.1 Assume $a_{13}=0,\,d_{13}\neq0$. From \eqref{A76} we get
$x_3=-x_2b_{23}^2/a_3$ and then since $a_{12}\neq0$ equation
\eqref{A102} yields $a_2=1$. From \eqref{A66} we get
$d_{13}=1-b_{23}^2c_{12}^2/a_3$. If we use \eqref{A81} to eliminate
second and higher powers of $a_3$ of equation \eqref{A82}, it factorizes
as $(a_3-1)(1+b_{23}c_{12})=0$, and we get two branches:

5.5.1.1
If we choose $a_3=1$ all other equations
are satisfied with $b_{23}=\omega^2/c_{12}$, where $\omega^3=1$ but we
must have $\omega\neq1$ to stay in the $(\beta,\gamma')$ case.
\beq \label{5.5.1.1}
    [1]  \begin{bmatrix}1-\omega& \frac{\omega }{ c}\\ c&.\end{bmatrix}
 \begin{bmatrix}.&.& \frac{\omega }{ c^{2}}
\\ x_{2}&1& \frac{-x_{2}  \omega }{ c^{2}}\\ c^{2}&.& 1-\omega 
 \end{bmatrix} \begin{bmatrix}.& \frac{\omega ^{2}}{ c}\\
   c  \omega ^{2}&1 -\omega\end{bmatrix} [1]
   \eeq
\medskip\noindent
{The eigenvalues are $1$ with multiplicity $6$ and $\omega$ with multiplicity $3$.}
\medskip

5.5.1.2 Now we choose $b_{23}=-1/c_{12}$ and then the remaining equations
are satisfied with $a_3=\varsigma=\pm i$.
\beq \label{eq:V.5.1.2}
    [1] \begin{bmatrix}\varsigma +1&
    \frac{-\varsigma }{ c}\\ c&.\end{bmatrix}
    \begin{bmatrix}.&.& \frac{-\varsigma }{ c^{2}}\\ x_{2}&1&
 \frac{x_{2}  \varsigma }{ c^{2}}\\ c^{2}&.&\varsigma  +1\end{bmatrix}
 \begin{bmatrix}.& \frac{-1}{ c}\\ \varsigma   c&\varsigma 
   +1\end{bmatrix} [\varsigma ]
   \eeq
   \medskip\noindent
   {The eigenvalues are $1$ with multiplicity $5$ and $\varsigma$ with multiplicity $4$.}
\medskip

5.5.2 The case $d_{13}=0,\,a_{13}\neq0$ is obtained by $02$-symmetry from 5.5.1.
{
Indeed, 
we see that the form $(\beta,\gamma^\prime)$ is invariant under the $\ket{0}\leftrightarrow\ket{2}$ symmetry, 
with the $3\times 3$ block having the following pairs interchanged $(a_{13},d_{13}),(b_{13},c_{13}),(x_2,x_3)$ and $(x_1,x_4)=(0,0)$.  Thus any solution obtained for $d_{13}=0$ and $a_{13}\neq 0$ may be transformed into a solution with $d_{13}\neq 0$ and $a_{13}=0$. 
}

\paragraph{Case 5.6} $(\beta,\delta')$

Since in this case $a_{12}\neq0$ we have from \eqref{A54} and
\eqref{A63} that $a_{13}=d_{13}=0$ but then \eqref{A78} implies
$a_{12}=0$, a contradiction.

\paragraph{Case 5.7} $(\delta,\delta')$.

From $\det\neq0$ we get $a_2\neq 0$ and then \eqref{A81} and
\eqref{A82} imply $c_{23}=c_{12}$ and $b_{23}=b_{12}$. Next from
\eqref{A38} and \eqref{A42} we get $c_{13}=c_{12}^2$ and
$b_{13}=b_{12}^2$. Equation \eqref{A39} then gives $a_3=1$ and
\eqref{A40} implies $c_{12}=1/b_{12}$. After this \eqref{A52} and
\eqref{A66} yield $a_{13}=d_{13}=0$. The remaining equations are satisfied
with $a_2=\epsilon,\,\epsilon=\pm 1$.
\beq\label{Case5.7} [1]  
  \begin{bmatrix}.&b\\ \frac{1}{ b}&.\end{bmatrix}
 \begin{bmatrix}.&.&b^{2}\\ x_{2}&\epsilon&x_{3}\\ \frac{1}{
 b^{2}}&.&.\end{bmatrix} \begin{bmatrix}.&b\\ \frac{1}{
 b}&.\end{bmatrix} [1]\eeq
\medskip\noindent
{ The eigenvalues are $1$ and $-1$ with multiplicities $5$ and $4$ if $\epsilon=-1$ and multiplicities $6$ and $3$ otherwise.}
\medskip

\subsection{Case 6: $x_1=x_2=x_3=0$, $x_4\neq0$ \label{ss:case6}}
From the outset it is best to divide this into two cases depending on whether or
not $b_{12}$ vanishes.

Case 6.1: $b_{12}=0$ therefore $a_{12}d_{12}\neq0$. Then from
\eqref{A1} $c_{12}=0$ and from \eqref{A29} and \eqref{A31}
$a_{12}=d_{12}=1$. From \eqref{A94} we get $b_{13}=0$, and hence
$a_{13}a_2d_{13}\neq0$ and then from \eqref{A90}, \eqref{A86} and
\eqref{A66} $a_{13}=a_2=d_{13}=1$, which leads to a contradiction with
\eqref{A68}.

Case 6.2: Now that $b_{12}\neq0$ we get from \eqref{A68} and
\eqref{A72} $a_{13}=d_{13}=0$. From \eqref{A72} $a_{13}=a_{12}$ and
from \eqref{A94} $b_{13}=b_{12}^2$ and then from \eqref{A98} and
\eqref{A99} $a_{12}=a_{23}=0$ and therefore $c_{12}b_{23}c_{23}\neq0$.
Next from \eqref{A46} $c_{13}=c_{12}^2$ and from \eqref{A100} and
\eqref{A101} $b_{23}=a_2^2/c_{12}$,  $b_{12}=a_2^2/c_{23}$ and from 
\eqref{A95} $a_3=c_{23}^2/c_{12}^2$.

At this point we divide the problem into two branches on whether or
not $d_{12}$ vanishes.

Case 6.2.1: $d_{12}=0$. Then from \eqref{A68} $d_{13}=0$ and from
\eqref{A95} $c_{23}=a_2^2c_{12}$ and from \eqref{A83}
$d_{23}=a_2(1-a_2^2)$. After this the remaining equations imply that
we must have either $a_2=\epsilon,$ $\epsilon=\pm 1$ or $a_2=\omega$
with $\omega^3=1$, $\omega\neq 1$. After changing $c_{12}\to c$ one
solution is
\beq \label{eq:VI.2.1a}
    [1] \begin{bmatrix}.& \frac{1}{
    c}\\ c&.\end{bmatrix}
 \begin{bmatrix}.&.& \frac{1}{ c^{2}}\\ .&\epsilon&.\\
   c^{2}&x_{4}&.\end{bmatrix}
   \begin{bmatrix}.& \frac{1}{
 c}\\ c& .
 \end{bmatrix} [1]
 \eeq
 \medskip\noindent
 {The eigenvalues are $1$ and $-1$, with multiplicities $5$ and $4$ if $\epsilon=-1$, otherwise multiplicities $6$ and $3$.  Note that this solution may be obtained from \eqref{Case5.7} by setting $x_2=0$ and taking the transpose, but this violates the Case 5 assumption that $x_2\neq 0$.}
 \medskip
 
 The second solution is
 \beq \label{eq:VI.2.1}
[1]  \begin{bmatrix}.& \frac{1}{ c}\\ c&.\end{bmatrix} 
 \begin{bmatrix}.&.& \frac{1}{ c^{2}}\\ .&\omega&.\\
   c^{2}&x_{4}&.\end{bmatrix}
   \begin{bmatrix}.& \frac{\omega^{2}}{
 c}\\ \omega^{2}  c&\omega-1
 \end{bmatrix} [ \omega]
 \eeq
\medskip\noindent
{The eigenvalues are $1,\omega$ and $-1$, each with multiplicity $3$.}
\medskip

Case 6.2.2: $d_{12}\neq0$. From \eqref{A64} and \eqref{A68} we get
$d_{13}=d_{23}=d_{12}(1-a_2)$ and then from \eqref{A63} $a_2=1$. The
remaining equations are solved by $d_{12}=(c_{23}-c_{12})/c_{23}$ and
$c_{23}=c_{12}\omega$ with $\omega^3=1$, $\omega\neq 1$, yielding:
\beq\label{sol:case14}
    [1]  \begin{bmatrix}.& \frac{1}{ c  \omega}\\ c& \omega+2\end{bmatrix}
      \begin{bmatrix}.&.& \frac{1}{ c^{2}  \omega^{2}}\\
 .&1&.\\ c^{2}&x_{4}&.\end{bmatrix}  \begin{bmatrix}.& \frac{1}{
 c}\\ \omega c&.\end{bmatrix}  [\omega^{2}]
\eeq
\medskip\noindent
{The eigenvalues are $1,\omega^2$ and $-\omega^2$, each with multiplicity $3$.}

\newcommand{\refx}{-}  \newcommand{\refxx}{}  \newcommand{\xf}[1]{\includegraphics[width=.2in]{xfig/#1}}
\newcommand{\evold}[2]{#1^{\times #2}}  \newcommand{\ev}[2]{{\begin{array}{c}#1\\{\times #2}\end{array}}}  \newcommand{\evo}[3]{$\left(\ev{1}{#1}\raisebox{.26cm}{,}\ev{{#2}}{#3}\right)$}
\newcommand{\evt}[5]{$\left(\ev{1}{#1}\raisebox{.27cm}{,}\ev{{#2}}{#3}\raisebox{.27cm}{,}\ev{{#4}}{#5}\right)$}
\newcommand{\rI}{1}
\newcommand{\rII}{2}
\newcommand{\rIII}{3}
\newcommand{\rIV}{4}
\newcommand{\rV}{5}
\newcommand{\rVI}{6}
\newcommand{\goo}[5]{#1 & #5 & #4 & #3 & #2 & \\ } 

\begin{table}[]
\begin{tabular}{l|lllll}
\goo{soln.   }{ eigenvalues, }{ parameters  }{ block         }{ non-zero  }   
\goo{name    }{ degeneracies }{ cont./discrete }{ form      }{ $x_i$s    }  
\hline
\goo{\rI     }{ \evo{8}{x}{1} }{ 3   /0  }{ \eqref{sol:case1} }{ \xf{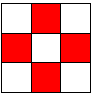}} 
\hline
\goo{\rII    }{      -      }{             }{  \refx{}      }{ \xf{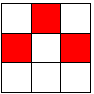}} 
\hline
\goo{\rIII.1.1}{\evt{5}{-x^2}{3}{x^3}{1}}{3 /0}{\eqref{eq:III.1.1}}{\xf{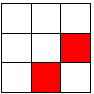} }
\goo{\rIII.1.2}{\evt{3}{-\omega x}{3}{\omega x^2}{3}}{3 /1}{\eqref{eq:III.1.2}}{ \xf{case3.pdf} }  
\hline
\goo{\rIV    }{      -    }{          }{       -         }{ \xf{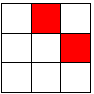}} 
\hline
\goo{\rV.2.1 }{\evo{7}{x}{2}}{3 /0}{\eqref{eq:V.2.1}   }{ \xf{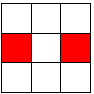} }
\goo{\rV.2.2 }{\evo{6}{x}{3}}{3 /0}{ \eqref{5.2.2}     }{ \xf{case5.pdf} }
\goo{\rV.4.a }{\evo{6}{x}{3}}{3 /0}{ \eqref{5.4.a}     }{\xf{case5.pdf}  }
\goo{\rV.4.b }{\evo{5}{x}{4}}{3 /0}{ \eqref{5.4.b}     }{ \xf{case5.pdf} }
\goo{\rV.5.1.1}{\evo{6}{\omega}{3}}{2 /1}{\eqref{5.5.1.1}}{ \xf{case5.pdf} }
\goo{\rV.5.1.2}{\evo{5}{\varsigma}{4}}{2 /1}{\eqref{eq:V.5.1.2}}{\xf{case5.pdf}}
\goo{\rV.7}{\evo{5}{-1}{4}/\evo{6}{-1}{3}}{3 /1}{\eqref{Case5.7}}{\xf{case5.pdf}} 
\hline
\goo{\rVI.2.1}{\evo{5}{-1}{4}/\evo{6}{-1}{3}}{2 /1}{\eqref{eq:VI.2.1a}}{ \xf{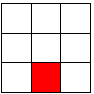} }
\goo{\rVI.2.1'}{\evt{3}{\omega}{3}{-1}{3}}{2 /1}{\eqref{eq:VI.2.1}}{\xf{case6.pdf} }
\goo{\rVI.2.2}{\evt{3}{\omega}{3}{-\omega}{3}}{2 /1}{\eqref{sol:case14}}{\xf{case6.pdf} }
\hline
\end{tabular}
\caption{Table of all ACC solutions
  to the Yang-Baxter equation (\ref{eq:braidr}) in rank-3. $\;$
  Here $x$ denotes a non-zero variable possibly with further constraints described in the text,
  $\omega$ is a primitive 3rd root of unity;
  and $\varsigma$ is a primitive 4th root of unity.
  In continuous/discrete parameter column entry
  3/1 means a 3-free-parameter family, not counting overall scaling, with 1 discrete parameter (which always take on exactly 2 values).
  (Hyphens and omitted `names' correspond to choices leading to no solution.)
  \label{tab:1}}
\end{table}

\medskip

Altogether we have established the following:

\begin{theorem}  \label{th:YB-ACC}
For the Yang-Baxter equation (\ref{eq:braidr}) in three dimensions,
the complete list of solutions satisfying ACC but not SCC
(see \cite{MR} for SCC)
is given, up to
noted symmetries (see Sec.\ref{ss:syms}),
in the formulae \eqref{sol:case1}-\eqref{sol:case14},
and collected in the Table~\ref{tab:1}.
\qed
\end{theorem}

\section{Analysis of the generic representations
  (\ref{eq:JaR1}/\ref{eq:boxx}).}
\label{ss:analysis}

Constant Yang--Baxter
solutions can be of considerable intrinsic interest. 
But they are also often interesting because of their symmetry algebras.
In the XXZ case 
(one of the strict charge-conserving cases)
for example, the symmetry
algebra is the quantum group $U_q sl_2$. 
This holds true in all ranks
(i.e. all system sizes $n$)  
- as we go up in ranks we 
simply see more of the symmetry algebra
-i.e. the action of the symmetry algebra on $n$-fold tensor space has a smaller kernel as $n$ increases.
It is not immediate that such a strong outcome would hold in general.
But it is interesting to investigate. 

In \S\ref{ss:analys} we  analyse our new solutions. 
(In \S\ref{ss:repthy1}
we recall
some classical facts about the classical cases for comparison.)

\subsection{Analysis of the generic solution: spectrum of $\R$}  \label{ss:analys}

\newcommand{\rr}{\check{r}}

Now we consider the solution in (\ref{eq:JaR1}/\ref{eq:boxx}). 
Observe that the trace of the $3\times3$ block is 
$$
a + \frac{\x_1 \x_3 +b}{b} -\frac{\x_1 \x_3}{ab} 
=
\frac{a^2 b + ab\x_1 \x_3 + ab^2 -\x_1 \x_3}{ab}
$$ 
Consider $\R_j -1$, so that all but the 3x3 block is zero. 
Restricting to the 3x3 block of $\R$, call it $\rr$, we have 
\beq  \label{eq:JaR12}
\rr -1_3 \; = \; 
\left[
\begin{array}{ccccccccc}
a & \x_1 & b
\\ 
 \frac{\mathit{\x_3} \left(a -1\right)}{b} 
 &  \frac{\mathit{\x_1} \mathit{\x_3} +b}{b} 
 &  \mathit{\x_3}
\\
\frac{\mathit{\x_3}^{2} \mathit{\x_1}^{2}}{b^{3}} 
& -\frac{\mathit{\x_1} \left(a b +\mathit{\x_1} \mathit{\x_3} \right)}{a \,b^{2}} 
& -\frac{\mathit{\x_3} \mathit{\x_1}}{a b} 
\end{array}\right] -1_3 \; 
= 
\left[ \begin{array}{ccccccccc}
a-1 & \x_1 & b
\\ 
 \frac{\mathit{\x_3} \left(a -1\right)}{b} 
 &  \frac{\mathit{\x_1} \mathit{\x_3} }{b} 
 &  \mathit{\x_3}
\\
\frac{\mathit{\x_3}^{2} \mathit{\x_1}^{2}}{b^{3}} 
& -\frac{\mathit{\x_1} \left(a b +\mathit{\x_1} \mathit{\x_3} \right)}{a \,b^{2}} 
& -\frac{b(\mathit{\x_3} \mathit{\x_1} +ab)}{a b^2}  
\end{array}\right]
\eeq 
Note that 
$\frac{ \x_3^2 \x_1^2 }{b^3} 
= \frac{-(a-1)(ab+\x_1 \x_3 )}{ab^2}
$
so this is clearly rank 1. Thus only one eigenvalue of $\R -1$ is not 0, and so only one eigenvalue of $\R$ is not 1.
We have 
\[
Trace(\R -1) = a-1 +\frac{\x_1 \x_3}{b} -\frac{b(\x_1 \x_3 +ab)}{ab^2}
= \frac{a^2 b^2 +ab\x_1 \x_3 -b(\x_1 \x_3 +ab)}{ab^2} -1 
\]
The other eigenvalue of $\R$ is 
$$
\lambda_2  = 
\frac{a(a-1) b^2 + (a-1)b \x_1 \x_3 }{ab^2}
= - (\x_1 \x_3 /b)^2
\; = -\left( \frac{-(a-1)\pm\sqrt{(a-1)^2-4a^2(a-1)}}{2a} \right)^2
$$
--- note from (\ref{eq:boxx}) that this depends only on $a$.

In particular each of our braid representations (varying the parameters appropriately) is a Hecke representation. 

\medskip

We see that the eigenvalue 
$\lambda_2$ can be varied over an open interval
(in each branch it is a continuous function of $a$, small for $a$ close to 1;
and large for large negative $a$). 
So (by Hecke representation theory,
specifically that the Hecke algebras are generically semisimple,
and abstract considerations
\cite{ClineParshallScott,Martin91})
the representation is generically semisimple. 

Returning to (\ref{eq:JaR12}) we have 
\[
\rr -1_3 \; =\;\; 
\left[ 1,\frac{\x_3}{b},\frac{-ab-\x_1 \x_3}{ab^2} \right]^t     
\;   [a-1,\x_1,b] 
=\;\; \frac{1}{ab^2}
\left[ 
\begin{array}{c}
ab^2 
\\ {ab\x_3} 
\\ {-ab-\x_1 \x_3} 
\end{array}
\right]     
\;   [a-1,\x_1,b] 
\]
and
\[
\R -1_9 \; =\;\; 
\left[0,0, 1,0,\frac{\x_3}{b},0,\frac{-ab-\x_1 \x_3}{ab^2},0,0 \right]^t     
\;   [0,0,a-1,0,\x_1,0,b,0,0] 
\]

Armed with this,  we have a Temperley--Lieb category representation
(i.e. an embedded TQFT - we assume familiarity with the standard $U_q sl_2$
version which can be used for comparison
- see e.g. \cite[Sec.6.2]{Bullivantetal2020} and references therein). 
In this form the duality is going to be
skewed (not a simple conjugation)
but should be workable.
In particular the loop parameter is
\[
[0,0,a-1,0,\x_1,0,b,0,0] 
\left[ 
\begin{array}{c}
0 \\ 0 \\ 1 \\ 0 
\\ {\x_3 /b} \\ 0  
\\ {(-ab-\x_1 \x_3)}/ab^2   \\ 0 \\ 0 
\end{array}
\right]  
=\;\; 
(a-1) + \frac{\x_1 \x_3}{b} -\frac{b(ab+\x_1 \x_3)}{ab^2}
\]\[ \hspace{2in} 
= \frac{ a^2 b^2  -2ab^2 +(a-1)b \x_1 \x_3}{ab^2}
= \frac{a-1}{a} (a + \frac{\x_1 \x_3}{b}) -1
= \lambda_2 -1
\]
which, note, depends only on $a$.

\subsection{Irreducible representation content of the generic solution $\rho_n$} 

The following analysis gives us
an invariant, and thus 
a way to classify solutions $\R$   (or equivalently $R$).

Thus in principle we can classify $R$-matrices according to the $B_n$-representation structure (the irrep content and so on) for each (and all) $n$. 
In general such an approach is very hard (due to the limits on knowledge of
the braid groups $B_n$ and their representation theory).
Certain properties can, however, make the problem more tractable. 

In our case call the representation $\rho_n$
(or just $\rho$ if no ambiguity arises,
or to denote the monoidal functor from the braid category, as in \cite{MR}). 
Depending on the field we are working over, this might mean the rep with indeterminate parameters, 
or a generic point in parameter space (i.e. the rep variety or a point on that variety).

Since this $\R$ has two eigenvalues (see \S\ref{ss:analys})
we have a Hecke representation
--- a representation of the algebra $H_n = H_n(q)$,
a quotient of the group algebra of $B_n$ for each $n$, for some $q$.
(With the same understanding about parameters.)

\newcommand{\aaaa}{\alpha}

Since eigenvalue $\lambda_1 =1$ 
this $H_n(q)$ is essentially in the `Lusztig' convention 
--- we can write $t_i$ for the braid generators in $H_n$, so 
$$
R_i = \rho(t_i) \;  = \;\; \rho_n (t_i)  \; ;
$$ 
then the quotient relation is 
\beq \label{eq:-1q} 
(t_i -1)(t_i +q)=0
\eeq  
for some $q=-\lambda_2$, as in \cite{PaulsNotes}. 
Here it is convenient to define
\[
U_i = \frac{t_i -1 }{ \aaaa}  
\]
so $\aaaa U_i (\aaaa U_i +1+q)=0$, 
i.e. 
$\aaaa U_i^2 = -(1+q) U_i$.

In a convention/parameterisation  as in (\ref{eq:-1q}) the operator 
\[
e' = 1-t_1 -t_2 +t_1 t_2 +t_2 t_1 -t_1 t_2 t_1
\]
is an unnormalised idempotent, and hence
\[
\rho_n(e') = 1-R_1 -R_2 +R_1 R_2 +R_2 R_1 -R_1 R_2 R_1
\]
is an unnormalised (possibly zero) idempotent,
whenever $\R$ gives such a Hecke representation. 

In our case in fact $e'$ is zero (by direct computation):
\beq  \label{eq:e00}
\rho_n(e' )  = 0
\eeq
Note that 
\[
\aaaa^3 U_1 U_2 U_1 = (t_1 -1)(t_2 -1)(t_1 -1) = 
t_1 t_2 t_1 
-t_1 t_2 -t_2 t_1 - t_1 t_1 
+2t_1 +t_2 -1
\]
so in our case 
\[
\rho_n( \aaaa^3 U_1 U_2 U_1 ) = (R_1 -1)(R_2 -1)(R_1 -1) = 
R_1 R_2 R_1 
-R_1 R_2 -R_2 R_1 - R_1 R_1 
+2R_1 +R_2 -1
\]
so \[
\rho_n( \aaaa^3 U_1 U_2 U_1 ) =  -R_1^2 +R_1 
= -R_1(R_1 -1)  
= \rho_n( q \aaaa  U_1 )
\hspace{1cm}\mbox{so}\hspace{1cm}
\rho_n( U_1 U_2 U_1 ) = \rho_n( \frac{q}{\aaaa^2} U_1 )
\]
so if we put $\aaaa = \pm \sqrt{q}$ then we have the relations of
the usual generators for Temperley--Lieb \cite{TemperleyLieb71}.

\medskip 

We assume familiarity with the generic irreducible representations of $H_n$,
which we write,
up to isomorphism,
as $L_\lambda$ with $\lambda\vdash n$ an integer partition of $n$.
The idempotent $e'$ induces the irrep $L_{1^3}$.
The unnormalised idempotent inducing the irrep  $L_3$ is
\beq \label{eq:en} 
e'_3 \; =\;  1+\frac{1}{q} (R_1 +R_2) +\frac{1}{q^2} (R_1 R_2 + R_2 R_1) +\frac{1}{q^3} R_1 R_2 R_1
\eeq 
This gives
\[
L_3(e'_3) = \frac{1}{q^3} (1+q)(1+q+q^2)
\]
which gives the normalisation factor, so 
\[
e_3 = \frac{q^3}{(1+q)(1+q+q^2)} e'_3
\]
The generalisation to irrep $L_n$ in rank $n$ will hopefully be clear
(in fact we won't really need it except for checking).

\medskip 

We can write $\chi_\lambda$ for the irreducible character associated to irrep $L_\lambda$.
That is, 
$$
\chi_{\lambda}(t_i) = Trace(L_{\lambda}(t_i)) .
$$
We can evaluate these characters in various ways, but a simple device is the restriction rule for the  inclusion $H_{n-1} \otimes 1_1 \hookrightarrow H_n$;
together with the easy cases:
\beq \label{eq:chin1n}
\chi_{n}(t_i)=1,  \hspace{1in} \chi_{1^n}(t_i)=\lambda_2
\eeq 
For example 
\[
\chi_{2,1}(t_i) = \chi_{2}(t_i)+\chi_{1^2}(t_i) =1+\lambda_2 
\]
and so on. 

Observe that the eigenvalues of $R_i$, specifically 
$R_1 = \R \otimes 1_3$, are three copies each of the eigenvalues of $\R$. 
Hence there are 24 eigenvalues $\lambda_1 =1$ and 3 copies of the other eigenvalue, call it $\lambda_2$:
\[
\chi_{\rho}(t_i) =3(8+\lambda_2)= 24+3\lambda_2
\]
The 1d irrep $L_3$, when present, contributes 1 eigenvalue $\lambda_1 =1$. 
The 2d irrep $L_{2,1}$ contributes 1 eigenvalue $\lambda_1 =1$ and 1 of the other eigenvalue $\lambda_2$.
The 1d irrep $L_{1^3}$ contributes just 1 of the other eigenvalue $\lambda_2$. 
Since $e'=0$ the multiplicity of this irrep in $\rho$ is 0. Therefore all the 3 eigenvalues $\lambda_2$ come from $L_{2,1}$ summands.
The identity (\ref{eq:e00}) therefore tells us that the irreducible content of our representation of $H_3$ (the Hecke quotient of $B_3$)
is 
\beq 
\rho  \; = \;\;  21 \; L_3 \; + \;\; 3\; L_{2,1} 
\eeq 
(the sum is generically but not necessarily always direct). 
In particular we have re-verified:

\noindent 
\propo{
Representation $\rho$ is a representation of Temperley--Lieb. 
}

\noindent

Note that it follows from the tensor construction that this TL property holds
(i.e. the image of $e'$ continues to vanish) for all $n$.
\\
Next we address the question of faithfulness of $\rho_n$ as a TL
representation,  and determine the
  centraliser, for all $n$.

\medskip 

Write $m_\lambda$ for the multiplicity of 
the generic irrep $L_\lambda$ in our rep $\rho$ (the generic character is well-defined in all specialisations, but the corresponding rep is not irreducible in all specialisations):
\beq \label{eq:mmultip}
\chi_{\rho_n} = \sum_{\lambda \vdash n} m_\lambda \chi_\lambda 
\eeq 

Note that integer partitions can be considered as vectors
(`weights' in Lie theory)
and hence added.
For example if $\mu=(\mu_1, \mu_2, \mu_3, ...,\mu_l)$ then
$$
\mu+11 = \mu+(1,1) = \mu+(1,1,0,...,0) = (\mu_1 +1, \mu_2 +1, \mu_3, ..., \mu_l ) .
$$

\medskip

\noindent 
{\bf{Stability Lemma}}. The multiplicity $m_\mu$ at level $n-2$ is the same as 
$m_{\mu+11}$ at level $n$. 

\medskip 

\noindent
{\em Outline Proof}.
The method of `virtual Lie theory' works here
(see e.g. \cite{Martin91,MartinRyomHansen}).
Let us define
\[
U_i = \R_i -1
\]
our rank-1 operator. Thus $U_i$ is itself an unnormalised idempotent - indeed it is, up to scalar, the image of the cup-cap operator in the TL diagram algebra. 

Write $T_n$ for TL on $n$ strands.
Recall that $U_1 T_n$ is a left $T_{n-2}$ right $T_n$ bimodule. Recall the algebra isomorphism $U_1 T_n U_1 \cong T_{n-2}$; and recall that 
$T_n / T_n U_1 T_n \cong k$ 
where $k$ is the ground field (for us it is $\C$). 
It follows that the category $T_{n-2}-mod$ embeds in 
$T_n-mod$, with embedding functor given by:
\beq \label{eq:TUM} 
M \mapsto T_n U_1 \otimes_{T_{n-2}} M
\eeq 
The irrep $L_\mu = L_{\mu_1, \mu_2}$ is taken to 
$L_{\mu+11} = L_{\mu_1 +1, \mu_2 +1}$.
Here $L_n$ is the module not hit by the embedding 
--- this is the module corresponding to  
$T_n / T_n U_1 T_n \cong k$, so the one that is annihilated by the localisation $M \mapsto U_1 M$. 
\hfill
\qed 

\medskip 

The Theorem below is a corollary of this Lemma.

It might also be of interest to show how to compute the further
multiplicities $m_\mu$ by direct calculation.
For $n=4$ we have 
$
\chi_{\rho_4}(t_i) = 3\chi_{\rho_3}(t_i) = 72+9\lambda_2 
.$
A direct calculation gives
$
\chi_{\rho_4}(e_4) = 55
$
so $m_4=55$, and we have
$
\chi_{\rho_4}(t_1) = 55 + m_{3,1} (2+\lambda_2) +m_{2,2} (1+\lambda_2).
$
We have $2m_{31} + m_{22} =72-55=17$ and 
$m_{31} + m_{22} =9$, giving $m_{3,1} = 8$ and $m_{2,2}=1$. 

Observe that this is in agreement with the Stability Lemma. 

For $n=5$ we have 
$
\chi_{\rho_5}(t_i) = 3\chi_{\rho_4}(t_i) = 216+27\lambda_2 
$
A direct calculation gives
$
\chi_{\rho_5}(e_5) = 144
$
and so we have
$
\chi_{\rho_5}(t_1) = 144 + m_{4,1} (3+\lambda_2) +m_{3,2} (3+2\lambda_2) 
$
We have $3m_{41} + 3m_{32} =216-144=72$ and 
$m_{41} + 2m_{32} =27$, giving $m_{4,1} = 21$ and $m_{3,2}=3$. 

We observe a pattern of repeated multiplicities, in agreement with the Stability Lemma: 
\[
\begin{array}{c|cccccccc}
m_\lambda & 1 & 3 & 8 & 21 & 55 & 144 \\ \hline 
  \lambda & 11&   & 2 \\
          &   &21 &   &3  \\
          & 22&   &31 &    & 4 \\
          &   &32 &   &41  &    & 5 \\
          & 33&   &   
\end{array}
\]
Besides the Stability Lemma or a direct calculation, 
the last entry above may be guessed based on Perron--Frobenius
applied to the Hamiltonian
$H = \sum_i \R_i$
--- if some power of $H$ is positive then there is a unique largest magnitude eigenvalue, and hence the corresponding multiplicity is 1. 
We know from the XXZ chain, which has the same eigenvalues but different multiplicities, that $\lambda = mm$ gives the largest eigenvalue when $n=2m$.

\begin{theorem}
The multiplicity $m_n$ 
in (\ref{eq:mmultip})
is given by A001906  from Sloane/OEIS \cite{Sloane},
with all other multiplicities $m_\mu$ determined by the Stability Lemma. 
\hfill\qed
\end{theorem} 

The
Temperley--Lieb algebras  are generically semisimple;
and a representation of a semisimple algebra is faithful
if and only if every irrep appears as a summand.
The latter is immediate from the Theorem,
so
generical faithfulness of our representations $\rho_n$ is similarly immediate. 

This brings us back to the original question about the stability of the
centraliser as $n$ varies - the possibility of an overarching symmetry algebra
analogous to $U_q sl_2$ in the XXZ case. Of course by Schur's Lemma the
Stability Lemma exactly says that there is a limit symmetry algebra, with all finite cases simply quotients of this limit.
But the combinatorial fact does not of itself imply that the symmetry algebra is something as beautiful as a quantum group (cf. Appendix~\ref{ss:repthy1}).

\appendix
\section{Appendix: The equations}
\subsection{The cubic constraints} \label{ss:cubics}

Here we write out the system of cubics corresponding to
entries in $\A$
as in (\ref{eq:AR}), 
hence the cubics that must vanish, in the ACC ansatz.

In fact the first few cubics in $\A$ are unchanged
(ordering 000 001 002 010 011 012 020 021 022 100 101 102 ... 222) from the strict CC ansatz.
Row 000 has vanishing anomaly.  
Row 001 gives:
\begin{equation*}
\bra{001} \A \ket{001} =\;
-\aaaaa_{12} \bbbbb_{12} \ccccc_{12}
-\aaaaa_1 \aaaaa_{12}^2
+\aaaaa_1^2 \aaaaa_{12},
\hspace{1.6cm} 
\bra{001}\A\ket{010}=
-\aaaaa_{12} \bbbbb_{12} \ddddd_{12}
\end{equation*}
with all other entries vanishing.
The first departure from SCC is in the 002 row,
which is:
\newcommand{\nougent}[1]{#1}  \[
\bra{002} \A \;\; = \;\; \; 
[0,0, \; 
-\aaaaa_{12} \xxxxx_1 \xxxxx_2
-\aaaaa_{13} \bbbbb_{13} \ccccc_{13}
-\aaaaa_1 \aaaaa_{13}^2
+\aaaaa_1^2 \aaaaa_{13},
\hspace{5.5cm} \]\[ 
\;\;0,\;\;\;
-\aaaaa_{13} \bbbbb_{13} \xxxxx_4
+(\aaaaa_1 \aaaaa_{12}
-\aaaaa_{12} \aaaaa_2
-\aaaaa_1 \aaaaa_{13} ) \xxxxx_1
,0, \; 
-\aaaaa_{12} \xxxxx_1 \xxxxx_3
-\aaaaa_{13} \bbbbb_{13} \ddddd_{13},
\;\; 0,0,0,
\]\[ 
\nougent{-\bbbbb_{13} \ccccc_{12} \xxxxx_1
-\aaaaa_{12} \bbbbb_{12} \xxxxx_1
+\aaaaa_1 \bbbbb_{12} \xxxxx_1}, \; 0,
\nougent{-\bbbbb_{13} \ddddd_{12} \xxxxx_1
+\aaaaa_1 \bbbbb_{13} \xxxxx_1
-\bbbbb_{12}^2 \xxxxx_1},
\;0,0,0,0,0,0,0,0,0,0,0,0,0,0]
\]


\allowdisplaybreaks

\numberwithin{equation}{section}
\setcounter{equation}{0}
\subsection{List of equations}
We give the complete list equations that are distinct up to an overall sign, organised by the number of terms (in computations we use the scale freedom to assume $a_1=1$).


\begin{eqnarray}
&&a_{12}\,c_{12}\,d_{12}=0,\label{A1}\\&&
a_{12}\,b_{12}\,d_{12}=0,\label{A2}\\&&
a_{23}\,c_{23}\,d_{23}=0,\label{A3}\\&&
a_{23}\,b_{23}\,d_{23}=0,\label{A4}\\&&
x_2\,x_4\,c_{12}=0,\label{A5}\\&&
x_2\,x_4\,c_{23}=0,\label{A6}\\&&
x_1\,x_3\,b_{12}=0,\label{A7}\\&&
x_1\,x_3\,b_{23}=0\label{A8}\\&&
a_{12}\,d_{12}\,(a_{12} - d_{12})=0,\label{A9}\\&&
a_{23}\,d_{23}\,(a_{23} - d_{23})=0,\label{A10}\\&&
x_1\,x_2\,(d_{12} - d_{23})=0,\label{A11}\\&&
x_1\,x_3\,(a_{12} - d_{12})=0,\label{A12}\\&&
x_1\,x_3\,(a_{23} - d_{23})=0,\label{A13}\\&&
x_2\,x_4\,(a_{12} - d_{12})=0,\label{A14}\\&&
x_2\,x_4\,(a_{23} - d_{23})=0,\label{A15}\\&&
x_3\,x_4\,(a_{12} - a_{23})=0,\label{A16}\\&&
x_1\,x_3\,c_{12} - a_{12}\,b_{12}\,d_{12}=0,\label{A17}\\&&
x_1\,x_3\,d_{23} + a_{13}\,b_{13}\,d_{13}=0,\label{A18}\\&&
x_1\,x_3\,a_{12} + a_{13}\,b_{13}\,d_{13}=0,\label{A19}\\&&
x_1\,x_3\,c_{23} - a_{23}\,b_{23}\,d_{23}=0,\label{A20}\\&&
x_1\,x_3\,d_{12} + a_{13}\,b_{13}\,d_{13}=0,\label{A21}\\&&
x_1\,x_3\,a_{23} + a_{13}\,b_{13}\,d_{13}=0,\label{A22}\\&&
x_2\,x_4\,d_{12} + a_{13}\,c_{13}\,d_{13}=0,\label{A23}\\&&
x_2\,x_4\,a_{12} + a_{13}\,c_{13}\,d_{13}=0,\label{A24}\\&&
x_2\,x_4\,b_{12} - a_{12}\,c_{12}\,d_{12}=0,\label{A25}\\&&
x_2\,x_4\,d_{23} + a_{13}\,c_{13}\,d_{13}=0,\label{A26}\\&&
x_2\,x_4\,a_{23} + a_{13}\,c_{13}\,d_{13}=0,\label{A27}\\&&
x_2\,x_4\,b_{23} - a_{23}\,c_{23}\,d_{23}=0,\label{A28}
\\&&
a_{12}\,(a_1^2 - a_1\,a_{12} - c_{12}\,b_{12})=0,\label{A29}\\&&
a_{23}\,(c_{23}\,b_{23} - a_3^2 + a_3\,a_{23})=0,\label{A30}\\&&
d_{12}\,(a_1^2 - a_1\,d_{12} - c_{12}\,b_{12})=0,\label{A31}\\&&
d_{23}\,(c_{23}\,b_{23} - a_3^2 + a_3\,d_{23})=0,\label{A32}\\&&
x_1\,(a_1\,b_{12} - c_{12}\,b_{13} - a_{12}\,b_{12})=0,\label{A33}\\&&
x_1\,(a_1\,b_{13} - d_{12}\,b_{13} - b_{12}^2)=0,\label{A34}\\&&
x_1\,(c_{23}\,b_{13} - a_3\,b_{23} + a_{23}\,b_{23})=0,\label{A35}\\&&
x_1\,(a_3\,b_{13} - d_{23}\,b_{13} - b_{23}^2)=0,\label{A36}\\&&
x_2\,(a_1\,c_{12} - c_{12}\,a_{12} - c_{13}\,b_{12})=0,\label{A37}\\&&
x_2\,(a_1\,c_{13} - c_{12}^2 - c_{13}\,d_{12})=0,\label{A38}\\&&
x_2\,(c_{13}\,a_3 - c_{13}\,d_{23} - c_{23}^2)=0,\label{A39}\\&&
x_2\,(c_{13}\,b_{23} - c_{23}\,a_3 + c_{23}\,a_{23})=0,\label{A40}\\&&
x_3\,(a_1\,b_{12} - c_{12}\,b_{13} - d_{12}\,b_{12})=0,\label{A41}\\&&
x_3\,(a_1\,b_{13} - a_{12}\,b_{13} - b_{12}^2)=0,\label{A42}\\&&
x_3\,(c_{23}\,b_{13} - a_3\,b_{23} + d_{23}\,b_{23})=0,\label{A43}\\&&
x_3\,(a_3\,b_{13} - a_{23}\,b_{13} - b_{23}^2)=0,\label{A44}\\&&
x_4\,(a_1\,c_{12} - c_{12}\,d_{12} - c_{13}\,b_{12})=0,\label{A45}\\&&
x_4\,(a_1\,c_{13} - c_{12}^2 - c_{13}\,a_{12})=0,\label{A46}\\&&
x_4\,(c_{13}\,a_3 - c_{13}\,a_{23} - c_{23}^2)=0,\label{A47}\\&&
x_4\,(c_{13}\,b_{23} - c_{23}\,a_3 + c_{23}\,d_{23})=0,\label{A48}
\\&&
x_3\,x_4\,(a_{12} - a_{23}) + x_1\,x_2\,( - d_{12} + d_{23})=0,\label{A49}\\&&
x_3\,x_4\,a_{23} - x_2\,x_1\,d_{23} + d_{13}\,a_{13}\,(d_{13} - a_{13})=0,\label{A50}\\&&
x_3\,x_4\,a_{12} - x_2\,x_1\,d_{12} + d_{13}\,a_{13}\,(d_{13} - a_{13})=0,\label{A51}\\&&
x_1\,x_2\,a_{12} + a_{13}\,( - a_1^2 + a_1\,a_{13} + c_{13}\,b_{13})=0,\label{A52}\\&&
x_1\,x_2\,a_{23} + a_{13}\,(c_{13}\,b_{13} - a_3^2 + a_3\,a_{13})=0,\label{A53}\\&&
x_1\,x_2\,b_{12} + b_{23}\,( - d_{23}\,a_{13} + a_{12}\,a_{13} - a_{12}\,a_{23})=0,\label{A54}\\&&
x_1\,x_2\,b_{23} + b_{12}\,( - d_{12}\,a_{13} - a_{12}\,a_{23} + a_{13}\,a_{23})=0,\label{A55}\\&&
x_1\,x_2\,c_{12} + c_{23}\,( - d_{23}\,a_{13} + a_{12}\,a_{13} - a_{12}\,a_{23})=0,\label{A56}\\&&
x_1\,x_2\,c_{23} + c_{12}\,( - d_{12}\,a_{13} - a_{12}\,a_{23} + a_{13}\,a_{23})=0,\label{A57}\\&&
x_1\,x_3\,a_2 + b_{13}\,(d_{13}\,a_{12} - d_{23}\,a_{12} + d_{23}\,a_{13})=0,\label{A58}\\&&
x_1\,x_3\,a_2 + b_{13}\,(d_{12}\,a_{13} - d_{12}\,a_{23} + d_{13}\,a_{23})=0,\label{A59}\\&&
x_2\,x_4\,a_2 + c_{13}\,(d_{12}\,a_{13} - d_{12}\,a_{23} + d_{13}\,a_{23})=0,\label{A60}\\&&
x_2\,x_4\,a_2 + c_{13}\,(d_{13}\,a_{12} - d_{23}\,a_{12} + d_{23}\,a_{13})=0,\label{A61}\\&&
x_3\,x_4\,b_{12} + b_{23}\,(d_{12}\,d_{13} - d_{12}\,d_{23} - d_{13}\,a_{23})=0,\label{A62}\\&&
x_3\,x_4\,b_{23} + b_{12}\,( - d_{12}\,d_{23} + d_{13}\,d_{23} - d_{13}\,a_{12})=0,\label{A63}\\&&
x_3\,x_4\,c_{12} + c_{23}\,(d_{12}\,d_{13} - d_{12}\,d_{23} - d_{13}\,a_{23})=0,\label{A64}\\&&
x_3\,x_4\,c_{23} + c_{12}\,( - d_{12}\,d_{23} + d_{13}\,d_{23} - d_{13}\,a_{12})=0,\label{A65}\\&&
x_3\,x_4\,d_{12} + d_{13}\,( - a_1^2 + a_1\,d_{13} + c_{13}\,b_{13})=0,\label{A66}\\&&
x_3\,x_4\,d_{23} + d_{13}\,(c_{13}\,b_{13} - a_3^2 + a_3\,d_{13})=0,\label{A67}\\&&
x_4\,(a_1\,d_{12} - a_1\,d_{13} - a_2\,d_{12}) - x_1\,c_{13}\,d_{13}=0,\label{A68}\\&&
x_4\,(a_2\,d_{23} + a_3\,d_{13} - a_3\,d_{23}) + x_1\,c_{13}\,d_{13}=0,\label{A69}\\&&
x_4\,a_{13}\,b_{13} + x_1\,(a_2\,a_{23} + a_3\,a_{13} - a_3\,a_{23})=0,\label{A70}\\&&
x_4\,a_{13}\,b_{13} + x_1\,( - a_1\,a_{12} + a_1\,a_{13} + a_2\,a_{12})=0,\label{A71}\\&&
x_4\,b_{12}\,(a_{12} - a_{13}) + x_1\,c_{12}\,( - d_{12} + d_{13})=0,\label{A72}\\&&
x_4\,b_{23}\,(a_{13} - a_{23}) + x_1\,c_{23}\,( - d_{13} + d_{23})=0,\label{A73}\\&&
x_2\,(a_1\,a_{12} - a_1\,a_{13} - a_2\,a_{12}) - x_3\,c_{13}\,a_{13}=0,\label{A74}\\&&
x_2\,(a_2\,a_{23} + a_3\,a_{13} - a_3\,a_{23}) + x_3\,c_{13}\,a_{13}=0,\label{A75}\\&&
x_2\,d_{13}\,b_{13} + x_3\,( - a_1\,d_{12} + a_1\,d_{13} + a_2\,d_{12})=0,\label{A76}\\&&
x_2\,d_{13}\,b_{13} + x_3\,(a_2\,d_{23} + a_3\,d_{13} - a_3\,d_{23})=0,\label{A77}\\&&
x_2\,b_{12}\,(d_{12} - d_{13}) + x_3\,c_{12}\,( - a_{12} + a_{13})=0,\label{A78}\\&&
x_2\,b_{23}\,(d_{13} - d_{23}) + x_3\,c_{23}\,( - a_{13} + a_{23})=0,\label{A79}\\&&
x_1\,(a_2\,b_{12} - a_2\,b_{23} + a_{12}\,b_{23} - a_{23}\,b_{12})=0,\label{A80}\\&&
x_2\,(c_{12}\,a_2 - c_{12}\,a_{23} - a_2\,c_{23} + c_{23}\,a_{12})=0,\label{A81}\\&&
x_3\,(a_2\,b_{12} - a_2\,b_{23} + d_{12}\,b_{23} - d_{23}\,b_{12})=0,\label{A82}\\&&
x_4\,(c_{12}\,a_2 - c_{12}\,d_{23} - a_2\,c_{23} + c_{23}\,d_{12})=0,\label{A83}\\&&
c_{12}\,d_{13}\,b_{12} - c_{23}\,d_{13}\,b_{23} + d_{12}^2\,d_{23} - d_{12}\,d_{23}^2=0,\label{A84}\\&&
c_{12}\,a_{13}\,b_{12} - c_{23}\,a_{13}\,b_{23} + a_{12}^2\,a_{23} - a_{12}\,a_{23}^2=0,\label{A85}\\
&&x_1\,x_2\,a_1 + x_3\,x_4\,a_{13} + a_{12}\,( - a_{12}\,a_2 - b_{12}\,c_{12} + a_2^2)=0,\label{A86}\\&&
x_1\,x_2\,a_3 + x_3\,x_4\,a_{13} + a_{23}\,( - a_{23}\,a_2 - b_{23}\,c_{23} + a_2^2)=0,\label{A87}\\&&
x_1\,x_2\,d_{13} + x_3\,x_4\,a_3 + d_{23}\,( - b_{23}\,c_{23} + a_2^2 - a_2\,d_{23})=0,\label{A88}\\&&
x_1\,x_2\,d_{13} + x_3\,x_4\,a_1 + d_{12}\,( - b_{12}\,c_{12} + a_2^2 - a_2\,d_{12})=0,\label{A89}\\&&
x_1\,x_2\,a_2 - c_{12}\,a_{23}\,b_{12} + c_{13}\,a_{23}\,b_{13} - d_{12}^2\,a_{13} + d_{12}\,a_{13}^2=0,\label{A90}\\&&
x_1\,x_2\,a_2 + c_{13}\,a_{12}\,b_{13} - c_{23}\,a_{12}\,b_{23} - d_{23}^2\,a_{13} + d_{23}\,a_{13}^2=0,\label{A91}\\&&
x_3\,x_4\,a_2 - c_{12}\,d_{23}\,b_{12} + c_{13}\,d_{23}\,b_{13} + d_{13}^2\,a_{12} - d_{13}\,a_{12}^2=0,\label{A92}\\&&
x_3\,x_4\,a_2 + c_{13}\,d_{12}\,b_{13} - c_{23}\,d_{12}\,b_{23} + d_{13}^2\,a_{23} - d_{13}\,a_{23}^2=0,\label{A93}\\&&
x_1\,(a_{13}\,d_{13} + a_2\,d_{12} - d_{12}\,d_{13}) + x_4\,( - b_{12}^2 + b_{13}\,a_1)=0,\label{A94}\\&&
x_1\,(a_{13}\,d_{13} + a_2\,d_{23} - d_{13}\,d_{23}) + x_4\,(b_{13}\,a_3 - b_{23}^2)=0,\label{A95}\\&&
x_1\,(a_1\,c_{13} - c_{12}^2) + x_4\,( - a_{12}\,a_{13} + a_{12}\,a_2 + a_{13}\,d_{13})=0,\label{A96}\\&&
x_1\,(c_{13}\,a_3 - c_{23}^2) + x_4\,( - a_{13}\,a_{23} + a_{13}\,d_{13} + a_{23}\,a_2)=0,\label{A97}\\&&
x_1\,(a_{13}\,d_{12} - b_{23}\,c_{12} + a_2^2 - a_2\,d_{12}) + x_4\,a_{23}\,b_{13}=0,\label{A98}\\&&
x_1\,(a_{13}\,d_{23} - b_{12}\,c_{23} + a_2^2 - a_2\,d_{23}) + x_4\,a_{12}\,b_{13}=0,\label{A99}\\&&
x_1\,c_{13}\,d_{12} + x_4\,( - a_{23}\,a_2 + a_{23}\,d_{13} - b_{23}\,c_{12} + a_2^2)=0,\label{A100}\\&&
x_1\,c_{13}\,d_{23} + x_4\,( - a_{12}\,a_2 + a_{12}\,d_{13} - b_{12}\,c_{23} + a_2^2)=0,\label{A101}\\&&
 x_2\,(a_1\,b_{13} - b_{12}^2) + x_3\,(a_2\,a_{12} + d_{13}\,a_{13} - a_{12}\,a_{13})=0,\label{A102}\\&&
x_2\,(a_3\,b_{13} - b_{23}^2) + x_3\,(a_2\,a_{23} + d_{13}\,a_{13} - a_{13}\,a_{23})=0,\label{A103}\\&&
x_2\,(a_2\,d_{12} - d_{12}\,d_{13} + d_{13}\,a_{13}) + x_3\,(a_1\,c_{13} - c_{12}^2)=0,\label{A104}\\&&
x_2\,(a_2\,d_{23} - d_{13}\,d_{23} + d_{13}\,a_{13}) + x_3\,(c_{13}\,a_3 - c_{23}^2)=0,\label{A105}\\&&
x_2\,d_{12}\,b_{13} + x_3\,(a_2^2 - a_2\,a_{23} - c_{23}\,b_{12} + d_{13}\,a_{23})=0,\label{A106}\\&&
x_2\,d_{23}\,b_{13} + x_3\,( - c_{12}\,b_{23} + a_2^2 - a_2\,a_{12} + d_{13}\,a_{12})=0,\label{A107}\\&&
x_2\,(c_{12}\,b_{23} - a_2^2 + a_2\,d_{23} - d_{23}\,a_{13}) - x_3\,c_{13}\,a_{12}=0,\label{A108}\\&&
x_2\,(a_2^2 - a_2\,d_{12} - c_{23}\,b_{12} + d_{12}\,a_{13}) + x_3\,c_{13}\,a_{23}=0,\label{A109}
 \end{eqnarray}

\section{Appendix:
Aside on further analysing solutions}
\label{ss:repthy1}

A step even further
than the all-ranks  representation theory analysis
in Sec.\ref{ss:analysis}
above would be to give an
{\em intrinsic}
characterisation of the centraliser algebra. We do not do this, but we can briefly set the scene. 

For an example $\R =P$ as in  (\ref{eq:brex})
is itself a solution --- this specific case, and also the corresponding $P$ for each $N$.
This solution is relatively simple, and completely understood in all cases, but still highly non-trivial. 
Of course it factors through the symmetric group.
(It is the Schur--Weyl dual to the natural general linear group action on tensor space.)
Its kernel as a symmetric group representation depends on $N$ as well as $n$. 
Assuming we work over the complex field, then the kernel is 
generated  exactly by the rank $N+1$ antisymmetriser.
Thus in particular for $N=2$ we have a faithful representation of `classical' Temperley--Lieb. While for $N=3$ the rank-3 antisymmetriser does not vanish
(so faithful on the corresponding algebras --- see e.g.
\cite{Brzezinski_1995}).

More explicitly we have the charge-conserving decomposition 
\[
\rho = (\rho_{111} \oplus \rho_{222}  \oplus \rho_{333} ) \oplus 
(\rho_{112}  \oplus \rho_{122}  \oplus \rho_{113}   \oplus \rho_{133}  \oplus \rho_{223}  \oplus \rho_{233} )
\oplus ( \rho_{123})
\]
\beq \label{eq:1081}
\cong 3\rho_{111} \oplus 6\rho_{112} \oplus \rho_{123}
\cong 10 L_3 \oplus 8 L_{21} \oplus L_{1^3}
\eeq 
where the bracketed sums are of isomorphic reps, and 
$\rho_{111}$ is trivial;
$\rho_{112} =L_3 \oplus L_{21}$;
$\rho_{123} = L_3 \oplus 2L_{21} \oplus L_{1^3}$
(i.e. the regular rep). 
Observe that the multiplicities 10, 8, 1 are the dimensions of the corresponding $GL_3$ irreps (recall these may be indexed by integer partitions of at most 2 rows, or equivalently of at most 3 rows where we delete all length-3 columns) as dictated by the duality.
Note that this structure will be preserved by any generic deformation. 

\newcommand{\beqa}{\begin{eqnarray}}
\newcommand{\eqa}{\end{eqnarray}} 
\newcommand{\yy}[1]{\begin{ytableau} #1 \end{ytableau}}
\newcommand{\ytwoone}{\yy{{}&{}\\{}}} 
\newcommand{\yoneone}{
\ytableausetup{smalltableaux} 
\begin{ytableau} { } \\ {} \end{ytableau} 
}

We can characterise this in the classical way, starting with the spectrum of $\R$ itself:
\begin{eqnarray}
\square  \otimes \square \; & = & \; \square\!\square \oplus \yoneone  \\ \label{eq:63}
3 \times 3 \;\; & = & \;\;\; 6 \;\;+\; \overline{3} 
\end{eqnarray}
\begin{eqnarray}
\square  \otimes \square\otimes\square \; & = & 
\; \left( \square\!\square \oplus \yoneone \right)\otimes\square   
\;\; = \;\; \yy{{}&{}&{}}\oplus 2. \yy{{}&{}\\{}} \oplus \emptyset 
\\ \label{eq:863}
3 \times 3\times 3 \;\; & = & 
\;\;\; (6 \;\;+\; \overline{3})\times 3 
\;\;\;\;\;\;\; = \;\; 10\;\; +\;\; 2.8 \;\;+ \;\;1
\end{eqnarray}
cf. (\ref{eq:1081}). 
Recall that this continues
\[
\square  \otimes \square\otimes\square\otimes\square  \;  = \yy{{}&{}&{}&{}} \oplus 3.\yy{{}&{}&{}\\{}} \oplus 2.\yy{{}&{}\\{}&{}} \oplus 3.\yy{{}}  
\;
\]
\[
3\times 3\times 3\times 3 \;\;= \;\; 15 \;\; +\;\; 3.15 \;\;+\;\; 2.\overline{6} \;\; +\;\; 3.3
\]
(Side note for future reference: Here in each third rank up the reps from three ranks down reappear (along with some more). This `three' is one sign that we are with $gl_3$ or $sl_3$ in this case.) 

\medskip 

Observe that
the  solution for $\R$ in (\ref{eq:JaR1}) 
(in \S\ref{ss:Jxxxx}) certainly does not have the multiplicities in (\ref{eq:63}). 
Indeed it agrees formally initially with
\begin{eqnarray*}
    \square \otimes \yoneone &=&  \ytwoone\oplus \emptyset 
    \\ 
    3\times \overline{3} \; &= &\;\;\; 8 \; + 1
\end{eqnarray*}
(see e.g. \cite{MartinRittenberg92})
- formally, in the sense that the symmetry needed for the symmetric group(/Hecke/braid) action is broken here. 
In this formal picture it is not clear how the labels would correspond with the Hecke algebra/symmetric group labels --- we are in rank-2 (but at least there are two summands). 
And it is not clear how to continue. We have 
\beqa 
\yy{{}} \otimes \ytwoone &=&  \yy{{}&{}&{}\\{}} \oplus \yy{{}&{}\\{}&{}} \oplus \yy{{}} 
\\
3 \times \; 8 \;\;\;&=&\;\;\; 15\;\; +\;\;\;  \overline{6} \;\;  + 3 
\eqa 
for example (so at least the centralised algebra of 
$ \yy{{}}\otimes\yy{{}\\{}}\otimes\yy{{}}$ 
is -miraculously - isomorphic to the Hecke quotient of $B_3$). 
But this is nowhere close to what we have.
This suggests that it is at least time to pass to the Lie supergroups again, such as $GL(2|1)$
(cf. e.g. \cite{Arnal94,Rittenberg,MartinRittenberg92}). 
(Alternatively it could be that the construction is not dual to a quantum group action.)

\bibliographystyle{abbrv}
\bibliography{local.bib,bib/Loop_Hecke.bib}

\end{document}